\documentclass[11pt]{article}
\usepackage[english]{babel}
\usepackage[cp1251]{inputenc}
\usepackage[psamsfonts]{amssymb}
\usepackage[intlimits]{amsmath}
\usepackage{amsfonts}
\usepackage{epsfig}
\usepackage{graphics}
\begin{document}

\begin{center}
{\bf ON A DIFFERENCE SCHEME FOR SOLVING CAUCHY PROBLEMS WITH THE CAPUTO FRACTIONAL DERIVATIVE IN A BANACH SPACE}\footnote{This work was supported by RFBR (project code 16--01--00039a)}
\\[2ex]
{\bf M.\,M.\,Kokurin}\\[3ex]
\end{center}

\medskip

We construct and study a time--semidiscretization scheme for the Cauchy problem associated with a linear homogeneous differential equation with the Caputo fractional time derivative of order $\alpha\in(0,1)$ and a spatial sectorial operator in a Banach space. For this scheme, we obtain rate--of--convergence and error estimates in terms of the discretization step. We use properties of Mittag--Leffler functions, hypergeometric functions, and the calculus of sectorial operators in a Banach space. Results of numerical experiments are also reported.

\medskip

\textbf{Keywords:} Cauchy problem, Caputo derivative, Banach space, finite difference scheme, error estimate, Mittag--Leffler function, hypergeometric function, operator calculus, sectorial operator.

\begin{center}
\textbf{1. Statement of the problem}
\end{center}

Let $u(t)$, $t\geq 0$ be a function with values in a complex
Banach space $X$. We use the standard
definition of the Riemann--Liouville fractional integral
$$
(I_{0+}^\alpha
u)(t)=\frac{1}{\Gamma(\alpha)}\int\limits_0^t\frac{u(s)}{(t-s)^{1-\alpha}}ds,\quad
t>0,
$$
and the Caputo fractional derivative $\partial^\alpha u=I_{0+}^{1-\alpha}u^{\prime}$ of order
$\alpha\in(0,1)$, see, e.g., [1,p.10--11; 2,Ch.2].

The paper is concerned with the Cauchy problem
$$
\partial^\alpha u(t)=Au(t),\quad t\in(0,T]; \quad u(0)=f\in D(A), \eqno{(1.1)}
$$
where $A: D(A) \subset X \to X$ is an unbounded closed operator and $\alpha\in(0,1)$. Following [3], by a solution of the problem (1.1) we mean a function $u\in C\bigl([0,T];X\bigr)$ that takes values in $D(A)$, satisfies the equations in (1.1) and such that $\partial^\alpha u(t)$ is a continuous function for $t>0$. The aim of this paper is to construct a finite difference method for approximate solution of (1.1) and to obtain estimates of the convergence rate of this method in terms of the discretization step.

Below we need the Mittag--Leffler function $E_\alpha(z)=\sum\limits_{n=0}^{\infty}\frac{z^n}{\Gamma(\alpha
n+1)}$, $\alpha>0$. The function $E_\alpha(z)=E_{\alpha,1}(z)$ is a special case of
the generalized Mittag--Leffler function
$E_{\alpha,\beta}(z)=\sum\limits_{n=0}^{\infty}\frac{z^n}{\Gamma(\alpha
n+\beta)}$, $\alpha>0$, $\beta\in\mathbb{C}$. In this paper, we will use $E_{\alpha,\beta}(z)$ with parameters $\alpha\in(0,1)$, $\beta\in\mathbb{R}$.

Several approaches to constructing difference schemes for scalar ($X=\mathbb{R}$) Cauchy problems with a fractional derivative are presented in [4; 5]. In [4; 6; 7], finite difference schemes for inhomogeneous fractional diffusion equations with specific spatial differential operators are constructed. In the present work, we suggest and study a scheme of semidiscretization in time for the problem (1.1) in a more general setting with an arbitrary sectorial operator $A$ in a Banach space, see Condition 1.1.

In [8], a finite difference scheme for the inhomogeneous Cauchy problem
$$
\partial^\alpha u(t)=Au(t)+g(t),\quad t\in(0,T], \quad u(0)=f
$$
is described. Here, $A$ is a differential operator of second order. The idea of constructing a difference scheme in the present paper and in [8] comes back to the finite difference approximation of the Caputo fractional derivative presented in [4]. However, the technique used in [4; 8] involves the assumption on $u(t)$ to be sufficiently smooth for $t\in[0,T]$. At the same time, even in the special case, where $A$ is the Laplace operator and $g(t)\equiv 0$, this condition is not satisfied, see Example 6.1. The reason is that the first time derivative of the function $E_\alpha(\lambda t^\alpha)$ is unbounded if $\lambda\neq 0$. In this paper, we suggest a way of overcoming the mentioned difficulty, see Section 2.

In [1,p.30], for the problem (1.1) with an operator $A$ of general form, a scheme of discrete approximation is described without accuracy estimates. An alternative approach to solving (1.1) based on a combination of the semidiscretization in time with a finite-dimensional approximation of spaces and operators, is developed in [9].

Below, $C_0$, $C_1$, $\dots$ are positive constants, which may depend on $\alpha$; $R(\zeta,T)=(\zeta E-T)^{-1}$ is the resolvent of the operator $T:D(T) \subset X\to X$, $\sigma(T)$ is its spectrum, $E$ is the identity operator in $X$. We also denote
$$
K(\varphi) = \{\zeta\in\mathbb{C}\backslash \{ 0\}\,|\, |\arg
\zeta|<\varphi\}.
$$
A linear operator $T:D(T)\subset X\to X$ is called
sectorial with the sectorial angle $\varphi_0\in(0,\pi)$ if
$\sigma(T)\subset K(\varphi_0)$ and
$$
\|R(\zeta,T)\|\leq\frac{C_0}{1+|\zeta|} \quad \forall \zeta
\in\mathbb{C}\backslash K(\varphi_0),
$$
where the constant $C_0$ does not depend on $\zeta$. Particularly, for a sectorial operator $T$ we have $T^{-1}\in L(X)$. Below we will use
the calculus of sectorial operators in a Banach space,
see, e.g., [10,Ch.2]. Suppose that a function $F(\zeta)$ is analytic in
$K(\varphi_0+\varepsilon)$ with $\varepsilon>0$ and
decreases at infinity, uniformly in $\arg\zeta$ for $\zeta\in K(\varphi_0+\varepsilon)$, faster than some negative power of
$|\zeta|$. Then the formula
$$
F(T) =\frac{1}{2\pi i}\int\limits_{\Gamma(r_0,\varphi_0)}
F(\zeta) R(\zeta,T)d\zeta, \eqno{(1.2)}
$$
with a sufficiently small $r_0>0$, defines the operator $F(T) \in
L(X)$. In this case,
$$
F(T)G(T) = (FG)(T) = \frac{1}{2\pi
i}\int\limits_{\Gamma(r_0,\varphi_0)} F(\zeta)G(\zeta)
R(\zeta, T) d\zeta. \eqno{(1.3)}
$$
The contour $\Gamma(r_0,\varphi_0)$ in (1.2), (1.3) with $\arg\zeta$ decreasing is defined as
\begin{equation}
\begin{split}
\Gamma(r_0,\varphi_0)=&\{\zeta=re^{-i\varphi_0}\,|\,r_0\leq
r<+\infty\}\cup\\
\cup&\{\zeta=r_0 e^{i\varphi}\,|\,|\varphi|\leq\varphi_0\}\cup\{\zeta=re^{i\varphi_0}\,|\,r_0\leq
r<+\infty\}.
\end{split}
\notag
\end{equation}
If $F(\zeta)=\frac{P_n(\zeta)}{Q_m(\zeta)}$, where $P_n(\zeta)$ and $Q_m(\zeta)$ are polynomials of degrees $n$ and $m$ respectively, $n<m$, and all zeroes of $Q_m(\zeta)$ lie outside of $\overline{K(\varphi_0)}$, then we have $F(T)=P_n(T)Q_m^{-1}(T)$. Throughout the paper, $\overline{D}$ denotes the closure of the set $D$.

The following condition specifies the class of problems (1.1) under consideration.

{\bf Condition 1.1.} {\it The operator $-A$ is sectorial with the sectorial angle $\varphi_0\in(0,\pi/2)$, and $\overline{D(A)}=X$.}

\medskip

If Condition 1.1 is satisfied, then the problem (1.1) is well--posed [3], in particular, there exists a unique solution $u(t)$, $t\in[0,T]$ for any $f\in D(A)$. In [11], we have established the following representation for the solution to (1.1) in terms of functions of the operator $-A$.

{\bf Lemma 1.1.} {\it The solution $u(t)$, $t\in(0,T]$ of (1.1) has the form $u(t)=F_t(-A)f$, where
$F_t(\zeta)=E_\alpha(-\zeta t^\alpha)$.}

\medskip

The paper is organized as follows. In Section 2 we describe a method of constructing difference approximations to the Caputo fractional derivative of a parametric family of scalar functions $E_\alpha(\lambda t^\alpha)$ on $[0,T]$. Using this approximation, in Section 3 we construct and study a finite difference scheme for solving a linear scalar Cauchy problem with the Caputo fractional derivative. In Section 4 we prove auxiliary integral inequalities used in Section 3. In Section 5 we construct a finite difference scheme for solving the Cauchy problem (1.1) and establish its accuracy estimates. The scheme is based on the scalar difference scheme introduced in Section 3, the construction also uses the calculus of sectorial operators. Section 6 provides results of numerical experiments. We compare the proposed solution method with a finite difference scheme from [1].

\begin{center}
\textbf{2. A finite difference approximation for the fractional Caputo derivative}
\end{center}

In this paper we follow the approach to constructing finite difference methods for the Cauchy problem in a Banach space, according to which at the first stage we study a difference scheme for the scalar analogue of the original problem. Then the obtained results are transferred to the corresponding difference scheme in a Banach space with the use of a suitable operator calculus. The outlined approach previously allowed to substantiate and investigate wide classes of stable finite difference methods for solving ill-posed Cauchy problems for operator differential equations of the first and second order in a Banach space [12; 13].

Let the operator $A$ in (1.1) satisfy Condition 1.1 with the sectorial angle $\varphi_0$. Consider the scalar analogue of (1.1), namely, the Cauchy problem for fractional ordinary--differential equations with the complex parameter $\lambda$:
$$
\partial^\alpha v(t)=\lambda v(t),\quad t\in(0,T];\quad v(0)=1. \eqno{(2.1)}
$$
In Sections 2 and 3 we assume that the following condition is fulfilled.

{\bf Condition 2.1.} {\it The parameter $\lambda$ in (2.1) belongs to the sector
$$
-\overline{K(\varphi_0)}=\{z\in\mathbb{C}\,|\, -z\in\overline{K(\varphi_0)}\}. \eqno{(2.2)}
$$}

The solution of (2.1) has the form [1,p.12; 3]
$$
v(t)=E_\alpha(\lambda t^\alpha). \eqno{(2.3)}
$$

For construction of difference schemes of solving (2.1), we need a discrete approximation for the fractional derivative $\partial^\alpha v(t_n)$ at the nodes of sampling $t_n=n\Delta t$, $0\leq n\leq N$ with $\Delta t=T/N$. The desired approximation is a linear combination of $v(t_0)$, $v(t_1)$, $\dots$, $v(t_n)$. For our purposes, it is sufficient to obtain such an approximation with a qualified error estimate for $v(t)$ of the form (2.3).

The classical approximations of fractional derivatives (see [4; 14]) are not applicable to $\partial^\alpha v(t)$ since $v(t)$ does not fulfill standard smoothness conditions, namely, the first derivative of $v(t)$ is unbounded on $[0,T]$.
The method [4] of constructing a difference approximation for $\partial^\alpha v(t)$ requires continuous differentiability of $v(t)$ on $[0,T]$ and is based on the following transformations:
$$
\partial^\alpha v(t_n)=\frac{1}{\Gamma(1-\alpha)}\int\limits_0^{t_n} \frac{v^\prime(s)ds}{(t_n-s)^\alpha} = \frac{1}{\Gamma(1-\alpha)}\sum\limits_{j=1}^n \int\limits_{t_{j-1}}^{t_j} \frac{v^\prime(s)ds}{(t_n-s)^\alpha} \approx
$$
$$
\approx \frac{1}{\Gamma(1-\alpha)}\sum\limits_{j=1}^n v^\prime(t_{j-1/2}) \int\limits_{t_{j-1}}^{t_j} \frac{ds}{(t_n-s)^\alpha} \approx \frac{1}{\Gamma(1-\alpha)}\sum\limits_{j=1}^n \frac{v(t_j)-v(t_{j-1})}{\Delta t} \int\limits_{t_{j-1}}^{t_j} \frac{ds}{(t_n-s)^\alpha}=
$$
$$
=\frac{1}{\Gamma(2-\alpha)}\sum\limits_{j=1}^n \Bigl((n-j+1)^{1-\alpha}-(n-j)^{1-\alpha}\Bigr) \frac{v(t_j)-v(t_{j-1})}{(\Delta t)^\alpha},
$$
$$
t_n=n\Delta t,\quad 0\leq n\leq N,\quad \Delta t=T/N.
$$
Below we modify this method using the identity
$$
v^\prime(t)=\frac{g(t)}{t^{1-\alpha}},\quad g(t)=\alpha\lambda E_\alpha^\prime(\lambda t^\alpha); \quad g\in C[0,T].
$$
Let us introduce additional nodes
$$
s_j\in(t_{j-1},t_j),\quad s_j^\alpha=\frac{t_{j-1}^\alpha+t_j^\alpha}{2},\quad 1\leq j\leq n.
$$
We have
$$
\partial^\alpha v(t_n)= \frac{1}{\Gamma(1-\alpha)}\sum\limits_{j=1}^n \int\limits_{t_{j-1}}^{t_j} \frac{g(s)ds}{s^{1-\alpha}(t_n-s)^\alpha} = \frac{1}{\Gamma(1-\alpha)}\sum\limits_{j=1}^n g(s_j)\int\limits_{t_{j-1}}^{t_j} \frac{ds}{s^{1-\alpha}(t_n-s)^\alpha} + \Delta_1 =
$$
$$
=\frac{1}{\Gamma(1-\alpha)}\sum\limits_{j=1}^n \bigl(\alpha\lambda E_\alpha^\prime(\lambda s_j^\alpha)\bigr) \int\limits_{t_{j-1}}^{t_j} \frac{ds}{s^{1-\alpha}(t_n-s)^\alpha}+\Delta_1=
$$
$$
=\frac{\alpha\lambda}{\Gamma(1-\alpha)}\sum\limits_{j=1}^n \Bigl( \frac{E_\alpha(\lambda t_j^\alpha)-E_\alpha(\lambda t_{j-1}^\alpha)}{\lambda t_j^\alpha-\lambda t_{j-1}^\alpha} + \Delta_2^j \Bigr) \int\limits_{j-1}^{j} \frac{dx}{x^{1-\alpha}(n-x)^\alpha}+\Delta_1=
$$
$$
=\frac{\alpha}{\Gamma(1-\alpha)}\sum\limits_{j=1}^n b_{jn} \frac{v(t_j)-v(t_{j-1})}{t_j^\alpha-t_{j-1}^\alpha} +\Delta_1+\Delta_2.
$$
Here, we have introduced the notation
$$
\Delta_1=\Delta_1(\lambda,n,\Delta t)=\frac{1}{\Gamma(1-\alpha)}\sum\limits_{j=1}^n \int\limits_{t_{j-1}}^{t_j}\frac{g(s)-g(s_j)}{s^{1-\alpha}(t_n-s)^\alpha}ds, \eqno{(2.4)}
$$
$$
\Delta_2^j=\Delta_2^j(\lambda,\Delta t)=E_\alpha^\prime(\lambda s_j^\alpha)-\frac{E_\alpha(\lambda t_j^\alpha)-E_\alpha(\lambda t_{j-1}^\alpha)}{\lambda t_j^\alpha-\lambda t_{j-1}^\alpha},\quad 1\leq j\leq n, \eqno{(2.5)}
$$
$$
\Delta_2=\Delta_2(\lambda,n,\Delta t)=\frac{\alpha\lambda}{\Gamma(1-\alpha)}\sum\limits_{j=1}^n b_{jn}\Delta_2^j(\lambda,\Delta t) \eqno{(2.6)}
$$
$$
b_{jn}=\int\limits_{j-1}^j \frac{dx}{x^{1-\alpha}(n-x)^\alpha},\quad 1\leq j\leq n\leq N.
$$
Thus, we obtain the approximation
$$
\partial^\alpha v(t_n) \approx \frac{\alpha}{\Gamma(1-\alpha)}\sum\limits_{j=1}^n b_{jn} \frac{v(t_j)-v(t_{j-1})}{t_j^\alpha-t_{j-1}^\alpha}. \eqno{(2.7)}
$$
It is easy to see that in the case of $\lambda=0$, (2.7) becomes the exact equality. Further assume that $\lambda\neq 0$. In this section, our goal is to establish estimates for the error $|\Delta_1+\Delta_2|$ of the approximation (2.7) in terms of $\lambda$ and $\Delta t$.

Let us estimate the coefficients $b_{jn}$. Obviously, they are all positive. We have
\begin{equation}
\begin{split}
&b_{jn}\leq\frac{1}{(j-1)^{1-\alpha}(n-j)^\alpha},\quad 2\leq j\leq n-1\\
&b_{1n}=\int\limits_0^1\frac{dx}{x^{1-\alpha}(n-x)^\alpha} \leq \frac{1}{(n-1)^\alpha}\int\limits_0^1\frac{dx}{x^{1-\alpha}} = \frac{1}{\alpha(n-1)^\alpha},\\
&b_{nn}=\int\limits_{n-1}^n\frac{dx}{x^{1-\alpha}(n-x)^\alpha} \leq \frac{1}{(n-1)^{1-\alpha}}\int\limits_{n-1}^n\frac{dx}{(n-x)^\alpha} = \frac{1}{(1-\alpha)(n-1)^{1-\alpha}}
\end{split}
\tag{2.8}
\end{equation}
when $2\leq n\leq N$ and
$$
b_{11}=\int\limits_0^1\frac{dx}{x^{1-\alpha}(1-x)^\alpha}=B(\alpha,1-\alpha).
$$

We proceed to obtain upper bounds for $|\Delta_1|$. We get
$$
|\Delta_1(\lambda,n,\Delta t)| \leq \frac{1}{\Gamma(1-\alpha)} \sum\limits_{j=1}^n\int\limits_{t_{j-1}}^{t_j}\frac{|g(s)-g(s_j)|}{s^{1-\alpha}(t_n-s)^\alpha}ds. \eqno{(2.9)}
$$
The following reasoning allows to estimate $|g(s)-g(s_j)|$ in (2.9). According to [15,Ch.3],
$$
\Bigl(\frac{d}{dt}\Bigr)^m E_\alpha(t^\alpha)=t^{-m}E_{\alpha,1-m}(t^\alpha),\quad m\geq 1. \eqno{(2.10)}
$$
It follows that
$$
v^{(m)}(s)=s^{-m}E_{\alpha,1-m}(\lambda, s^\alpha),\quad m\geq 1. \eqno{(2.11)}
$$
Thus,
$$
g(s)=s^{1-\alpha}v^\prime(s)=G(s^\alpha)
$$
with the function $G(t)=t^{-1}E_{\alpha,0}(\lambda t)$ which is analytic in $t$. Let us find the derivative of $G(t)$. Using the substitution $\tau=t^{1/\alpha}$ and (2.11), we obtain
$$
G^\prime(t) = \frac{d}{dt}\bigl(t^{-1}E_{\alpha,0}(\lambda t)\bigr) = \frac{d}{d\tau}\bigl(\tau^{-\alpha}E_{\alpha,0}(\lambda\tau^\alpha)\bigr)\cdot\frac{1}{\alpha}t^{\frac{1}{\alpha}-1} = \frac{d}{d\tau}\bigl(\tau^{1-\alpha}v^\prime(\tau)\bigr)\cdot\frac{1}{\alpha}t^{\frac{1}{\alpha}-1}=
$$
$$
= \frac{1}{\alpha}\tau^{-\alpha} \bigl( (1-\alpha)v^\prime(\tau)+\tau v^{\prime\prime}(\tau) \bigr)t^{\frac{1}{\alpha}-1} = \frac{1}{\alpha}t^{-2} \bigl( (1-\alpha)E_{\alpha,0}(\lambda t)+E_{\alpha -1}(\lambda t) \bigr).
$$

We claim that for all $\lambda$ satisfying Condition 2.1,
$$
|G^\prime(t)|\leq C_1 |\lambda|^2, \eqno{(2.12)}
$$
with $C_1$ independent of $t$ and $\lambda$. In fact,
$$
\frac{G^\prime(t)}{\lambda^2}=\widetilde G(\lambda t),\quad \widetilde G(z) = \frac{1}{\alpha z^2} \bigl( (1-\alpha)E_{\alpha,0}(z)+E_{\alpha -1}(z) \bigr).
$$
We have the series expansions:
$$
v(s) = E_\alpha(\lambda s^\alpha) = \sum\limits_{n=0}^\infty \frac{\lambda^n s^{\alpha n}}{\Gamma(\alpha n+1)},
$$
$$
g(s) = s^{1-\alpha}v^\prime(s) = \sum\limits_{n=0}^\infty \frac{\alpha(n+1) \lambda^{n+1} s^{\alpha n}}{\Gamma(\alpha n+\alpha+1)} = G(s^\alpha),
$$
$$
G(t) = \sum\limits_{n=0}^\infty \frac{\alpha(n+1)\lambda^{n+1}t^n}{\Gamma(\alpha n+\alpha+1)},\quad G^\prime(t) = \sum\limits_{n=0}^\infty \frac{\alpha(n+1)(n+2)\lambda^{n+2}t^n}{\Gamma(\alpha n +2\alpha+1)} = \lambda^2\widetilde G(\lambda t),
$$
$$
\widetilde G(z) = \sum\limits_{n=0}^\infty \frac{\alpha(n+1)(n+2)z^n}{\Gamma(\alpha n +2\alpha+1)}.
$$
We see that $\widetilde G(z)$ is an entire function. In accordance with Condition 2.1, we consider $\widetilde G(z)$ on $-\overline{K(\varphi_0)}$. By the asymptotic expansions from [11; 16,Ch.1], it follows that
$$
|E_{\alpha,\beta}(z)|\leq C_2 |z|^{-1} \eqno{(2.13)}
$$
for arbitrary $\alpha\in(0,1)$, $\beta\in\mathbb{R}$ and $z$ in the sector $|\arg z|>\alpha\pi/2+\varepsilon$ with sufficiently large $|z|$. Consequently the function $\widetilde G(z)$ is bounded on this sector, and hence on $-\overline{K(\varphi_0)}$. This yields the desired estimate (2.12).

Let us turn to evaluating the expression $|g(s)-g(s_j)|$, $s\in[t_{j-1},t_j]$ in (2.9). By the Lagrange theorem [17,Ch.1; 18,Ch.8],
$$
|g(s)-g(s_j)|=|G(s^\alpha)-G(s_j^\alpha)| \leq \max\limits_{t\in [s^\alpha, s_j^\alpha]}|G^\prime(t)| \cdot |s^\alpha-s_j^\alpha|.
$$
Combining the last inequality with (2.12), we obtain
$$
|g(s)-g(s_j)| \leq C_1|\lambda|^2 |s^\alpha - s_j^\alpha| \leq C_1|\lambda|^2 \frac{t_j^\alpha-t_{j-1}^\alpha}{2} = C_1|\lambda|^2 (\Delta t)^\alpha \frac{j^\alpha-(j-1)^\alpha}{2}.
$$
It is easily seen that
$$
j^\alpha-(j-1)^\alpha\leq\frac{C_3}{j^{1-\alpha}},\quad 1\leq j\leq n. \eqno{(2.14)}
$$
Hence we come to the estimate
$$
|g(s)-g(s_j)|\leq\frac{C_4|\lambda|^2(\Delta t)^\alpha}{j^{1-\alpha}}. \eqno{(2.15)}
$$

By substituting (2.15) into (2.9) and using the definition of coefficients $b_{jn}$, $1\leq j\leq n\leq N$, we get
$$
|\Delta_1(\lambda,n,\Delta t)| \leq C_5|\lambda|^2(\Delta t)^\alpha \sum\limits_{j=1}^n \frac{b_{jn}}{j^{1-\alpha}}. \eqno{(2.16)}
$$
In view of (2.8), the sum in (2.16) without the first two and last two terms is estimated as follows
$$
\sum\limits_{j=3}^{n-2}\frac{b_{jn}}{j^{1-\alpha}} \leq C_6\sum\limits_{j=3}^{n-2}f(j);\quad f(x)=(x-1)^{-2+2\alpha}(n-x)^{-\alpha},\quad x\in(1,n). \eqno{(2.17)}
$$
The function $f(x)$ goes to $+\infty$ when $x\to 1$, $x\to n$ and has a single minimum on the interval $(1,n)$. We now estimate the sum in (2.17) by the corresponding integral and then make the change of variable $x=(n-1)t+1$:
$$
\sum\limits_{j=3}^{n-2}\frac{b_{jn}}{j^{1-\alpha}} \leq C_6\int\limits_2^{n-1}f(x)dx = \frac{C_6}{(n-1)^{1-\alpha}}\int\limits_{\frac{1}{n-1}}^{\frac{n-2}{n-1}}\frac{dt}{t^{2-2\alpha}(1-t)^\alpha}\leq
$$
$$
\leq\frac{C_6}{(n-1)^{1-\alpha}} \biggl( \int\limits_{\frac{1}{n-1}}^{1/2}\frac{dt}{t^{2-2\alpha}(1-t)^\alpha} + \int\limits_{1/2}^1\frac{dt}{t^{2-2\alpha}(1-t)^\alpha} \biggr)\leq \frac{C_6}{(n-1)^{1-\alpha}} \biggl( C_7\int\limits_{\frac{1}{n-1}}^{1/2}\frac{dt}{t^{2-2\alpha}} + C_8\biggr).
$$
Estimating the omitted first two and last two summands, we come to the inequality
\begin{equation}
\sum\limits_{j=1}^n\frac{b_{jn}}{j^{1-\alpha}} \leq C_9
\begin{Bmatrix}
n^{-\alpha},&\alpha\in(0,1/2);\\
n^{-1/2}\ln(n+1),&\alpha=1/2;\\
n^{-1+\alpha},&\alpha\in(1/2,1)
\end{Bmatrix}
=C_9\eta(n).
\tag{2.18}
\end{equation}
It is easy to see that (2.18) is true for all $1\leq n\leq N$, including those $n$ for which the sum (2.17) is not defined. Substituting (2.18) into (2.16), we arrive at the following assertion.

{\bf Lemma 2.1.} {\it For the value $\Delta_1(\lambda,n,\Delta t)$ defined in (2.4), under Condition 2.1 we have
$$
|\Delta_1(\lambda,n,\Delta t)|\leq C_{10}\eta(n)|\lambda|^2(\Delta t)^\alpha.
$$}

\medskip

Let us obtain estimates for $|\Delta_2^j|$, $1\leq j\leq n$ from (2.5). Applying the Taylor formula with the integral remainder term [17,Ch.1; 19] and using the equality $t_j^\alpha s_j^\alpha=s_j^\alpha-t_{j-1}^\alpha=(t_j^\alpha-t_{j-1}^\alpha)/2$, we can express $E_\alpha(\lambda t_j^\alpha)$ and $E_\alpha(\lambda t_{j-1}^\alpha)$ by the values of the Mittag--Leffler function and its derivatives at the point $\lambda s_j^\alpha$ as follows:
$$
E_\alpha(\lambda t_j^\alpha) = E_\alpha(\lambda s_j^\alpha) + \frac{E_\alpha^\prime(\lambda s_j^\alpha)}{2} \lambda(t_j^\alpha-t_{j-1}^\alpha) +
$$
$$
+ \frac{E_\alpha^{\prime\prime}(\lambda s_j^\alpha)}{8} \lambda^2(t_j^\alpha-t_{j-1}^\alpha)^2 + \frac{1}{2} \int\limits_{[\lambda s_j^\alpha, \lambda t_j^\alpha]} E_\alpha^{\prime\prime\prime}(\tau)(\lambda t_j^\alpha-\tau)^2 d\tau;
$$
$$
E_\alpha(\lambda t_{j-1}^\alpha) = E_\alpha(\lambda s_j^\alpha) - \frac{E_\alpha^\prime(\lambda s_j^\alpha)}{2} \lambda(t_j^\alpha-t_{j-1}^\alpha) +
$$
$$
+\frac{E_\alpha^{\prime\prime}(\lambda s_j^\alpha)}{8} \lambda^2(t_j^\alpha-t_{j-1}^\alpha)^2 + \frac{1}{2} \int\limits_{[\lambda s_j^\alpha, \lambda t_{j-1}^\alpha]} E_\alpha^{\prime\prime\prime}(\tau)(\tau-\lambda t_{j-1}^\alpha)^2 d\tau.
$$
Consequently,
\begin{equation}
\begin{split}
|\Delta_2^j(\lambda,\Delta t)| = \Bigl| E_\alpha^\prime(\lambda s_j^\alpha)-&\frac{E_\alpha(\lambda t_j^\alpha)-E_\alpha(\lambda t_{j-1}^\alpha)}{\lambda t_j^\alpha-\lambda t_{j-1}^\alpha} \Bigr|=\\
=\frac{1}{2|\lambda|(t_j^\alpha-t_{j-1}^\alpha)} \biggl| \int\limits_{[\lambda s_j^\alpha, \lambda t_{j-1}^\alpha]} E_\alpha^{\prime\prime\prime}(\tau) (\tau-\lambda t_{j-1}^\alpha)^2 d\tau - &\int\limits_{[\lambda s_j^\alpha, \lambda t_j^\alpha]} E_\alpha^{\prime\prime\prime}(\tau) (\lambda t_j^\alpha-\tau)^2 d\tau \biggr| \leq\\
\leq\frac{1}{2|\lambda|(j^\alpha-(j-1)^\alpha)(\Delta t)^\alpha} \biggl( &\int\limits_{[\lambda s_j^\alpha, \lambda t_{j-1}^\alpha]} |E_\alpha^{\prime\prime\prime}(\tau)| \cdot |\tau-\lambda t_{j-1}^\alpha|^2 |d\tau|+\\
+&\int\limits_{[\lambda s_j^\alpha, \lambda t_j^\alpha]} |E_\alpha^{\prime\prime\prime}(\tau)|\cdot |\lambda t_j^\alpha-\tau|^2 |d\tau| \biggr).
\end{split}
\tag{2.19}
\end{equation}

Using (2.10) it is easy to get
$$
E_\alpha^{\prime\prime\prime}(\tau) = \frac{1}{\alpha^3}\tau^{-3} \bigl( E_{\alpha -2}(\tau)+3(1-\alpha)E_{\alpha -1}(\tau)+(1-\alpha)(1-2\alpha)E_{\alpha,0}(\tau) \bigr).
$$
Obviously, $E_\alpha^{\prime\prime\prime}(\tau)$ is an entire function. Therefore, by (2.13) it follows that $E_\alpha^{\prime\prime\prime}(\tau)$ is bounded in the sector (2.2). From (2.19) by direct calculations with the use of (2.14), we obtain
\begin{equation}
\begin{split}
&|\Delta_2^j(\lambda,\Delta t)| \leq \frac{C_{11}}{(j^\alpha-(j-1)^\alpha)|\lambda|(\Delta t)^\alpha} \biggl( \int\limits_{[\lambda s_j^\alpha, \lambda t_{j-1}^\alpha]} |\tau-\lambda t_{j-1}^\alpha|^2 |d\tau| +\\
&+\int\limits_{[\lambda s_j^\alpha, \lambda t_j^\alpha]} |\lambda t_j^\alpha-\tau|^2 |d\tau| \biggr)=\frac{C_{11}}{12} \bigl( j^\alpha-(j-1)^\alpha \bigr)^2 |\lambda|^2 (\Delta t)^{2\alpha} \leq \frac{C_{12}|\lambda|^2 (\Delta t)^{2\alpha}}{j^{2-2\alpha}}.
\end{split}
\tag{2.20}
\end{equation}

Now we can estimate $|\Delta_2|$. Using (2.8) and (2.20), we derive
$$
|\Delta_2(\lambda,n,\Delta t)|= \biggl|\frac{\alpha\lambda}{\Gamma(1-\alpha)}\sum\limits_{j=1}^n b_{jn}\Delta_2^j(\lambda,\Delta t)\biggr| \leq C_{13}|\lambda|\sum\limits_{j=1}^n b_{jn} |\Delta_2^j(\lambda,\Delta t)| \leq
$$
$$
\leq C_{14}|\lambda|^3 (\Delta t)^{2\alpha} \biggl( \frac{1}{n^\alpha} +\frac{1}{n^{3-3\alpha}} + \sum\limits_{j=2}^{n-1}\frac{1}{(j-1)^{3-3\alpha}(n-j)^\alpha} \biggr).
$$
Here, the sum in brackets is estimated similarly to (2.17):
$$
\sum\limits_{j=2}^{n-1}\frac{1}{(j-1)^{3-3\alpha}(n-j)^\alpha} \leq \frac{1}{(n-2)^\alpha}+ \frac{1}{(n-2)^{3-3\alpha}}+ \int\limits_2^{n-1}\frac{dx}{(x-1)^{3-3\alpha}(n-x)^\alpha}.
$$
To estimate the integral, we make the substitution $x=(n-1)t+1$:
$$
\int\limits_2^{n-1}\frac{dx}{(x-1)^{3-3\alpha}(n-x)^\alpha}= \frac{1}{(n-1)^{2-2\alpha}}\int\limits_\frac{1}{n-1}^\frac{n-2}{n-1}\frac{dt}{t^{3-3\alpha}(1-t)^\alpha}\leq
$$
$$
\leq \frac{1}{(n-1)^{2-2\alpha}} \biggl( \int\limits_\frac{1}{n-1}^{1/2}\frac{dt}{t^{3-3\alpha}(1-t)^\alpha} + \int\limits_{1/2}^1\frac{dt}{t^{3-3\alpha}(1-t)^\alpha} \biggr) \leq
$$
\begin{equation}
\leq \frac{1}{(n-1)^{2-2\alpha}} \biggl( C_{15} \int\limits_\frac{1}{n-1}^{1/2}\frac{dt}{t^{3-3\alpha}} + C_{16} \biggr) \leq C_{17}
\begin{Bmatrix}
n^{-\alpha},&\alpha\in(0,2/3);\\
n^{-2/3}\ln(n+1),&\alpha=2/3;\\
n^{-2+2\alpha},&\alpha\in(2/3,1)
\end{Bmatrix}
=C_{17}\widetilde\eta(n).
\notag
\end{equation}
Thus,
$$
|\Delta_2(\lambda,n,\Delta t)|\leq C_{14}|\lambda|^3 (\Delta t)^{2\alpha} \Bigl( \frac{1}{n^\alpha} +\frac{1}{n^{3-3\alpha}} + \frac{1}{(n-2)^\alpha}+ \frac{1}{(n-2)^{3-3\alpha}}+ C_{17}\widetilde\eta(n) \Bigr) \leq
$$
$$
\leq C_{18}\widetilde\eta(n)|\lambda|^3(\Delta t)^{2\alpha}.
$$
It is easy to see that the obtained estimate is true for all $1\leq n\leq N$, although some intermediate expressions above have no sense for $n=1,2$. We have proved the following lemma.

{\bf Lemma 2.2.} {\it For the value $\Delta_2(\lambda,n,\Delta t)$ defined in (2.6), under Condition 2.1 we have
$$
|\Delta_2(\lambda,n,\Delta t)| \leq C_{18}\widetilde\eta(n)|\lambda|^3(\Delta t)^{2\alpha}.
$$
}

\medskip

Combining Lemmas 2.1 and 2.2 we arrive at the following result.

{\bf Lemma 2.3.} {\it Let Condition 2.1 be fulfilled. Then for the Caputo fractional derivative of the function (2.3), we have the estimate:
$$
\biggl|\partial^\alpha v(t_n)-\frac{\alpha}{\Gamma(1-\alpha)}\sum\limits_{j=1}^n b_{jn} \frac{v(t_j)-v(t_{j-1})}{t_j^\alpha-t_{j-1}^\alpha}\biggr|\leq C_{19}\eta(n)|\lambda|^2(\Delta t)^{\alpha}(1+|\lambda|(\Delta t)^\alpha).
$$
Here, $\eta(n)$ is defined in (2.18).
}

\medskip

In the next section, the approximation (2.7) will be used for constructing a finite difference scheme for the scalar Cauchy problem (2.1).

\begin{center}
\textbf{3. A difference scheme for the scalar Cauchy problem with\\ the Caputo fractional derivative}
\end{center}

We write the difference approximation (2.7) for the fractional Caputo derivative of the function (2.3) as
$$
\partial^\alpha v(t_n) \approx \sum\limits_{j=0}^n \frac{a_{jn}v(t_j)}{(\Delta t)^\alpha}.
$$
Here,
\begin{equation}
\begin{split}
a_{0n}&=-\frac{\alpha}{\Gamma(1-\alpha)}b_{1n},\\
a_{jn}&=\frac{\alpha}{\Gamma(1-\alpha)}\Bigl( \frac{b_{jn}}{j^\alpha-(j-1)^\alpha} - \frac{b_{j+1,n}}{(j+1)^\alpha j^\alpha} \Bigr),\quad 1\leq j\leq n-1 \\
a_{nn}&=\frac{\alpha}{\Gamma(1-\alpha)}\frac{b_{nn}}{n^\alpha-(n-1)^\alpha};\quad\quad 1\leq n\leq N.
\end{split}
\tag{3.1}
\end{equation}
Using this approximation, we arrive at the following finite difference scheme for the parametric Cauchy problem (2.1):
$$
\sum\limits_{j=0}^n \frac{a_{jn}v_j}{(\Delta t)^\alpha}-\lambda v_n=0,\quad 1\leq n\leq N;\quad v_0=1. \eqno{(3.2)}
$$
In (3.2), $v_n$ is the desired approximation to $v(t_n)$. In this section, one of our goals is to obtain rate--of--convergence estimates for the scheme (3.2) in terms of $\lambda$ and $\Delta t$. We emphasize that the scheme (3.2) is designed for scalar fractional derivative problems of the highly specialized class (2.1). However the results on the rate of convergence of (3.2) will be used in Section 5 for construction and justification of difference method for solving the Cauchy problem (1.1) with an arbitrary sectorial operator $A$ in a Banach space.

It is clear that $a_{0n}<0$ and $a_{nn}>0$ for all $1\leq n\leq N$. To determine the sign of remaining coefficients $a_{jn}$, we write the chain of inequalities:
$$
\frac{1}{(n-j+1)^\alpha}\int\limits_{j-1}^j \frac{dx}{x^{1-\alpha}} \leq \int\limits_{j-1}^j \frac{dx}{x^{1-\alpha}(n-x)^\alpha} \leq \frac{1}{(n-j)^\alpha}\int\limits_{j-1}^j \frac{dx}{x^{1-\alpha}},\quad 1\leq j\leq n;
$$
$$
\frac{1}{(n-j+1)^\alpha}\frac{j^\alpha-(j-1)^\alpha}{\alpha}\leq b_{jn}\leq \frac{1}{(n-j)^\alpha}\frac{j^\alpha-(j-1)^\alpha}{\alpha},\quad 1\leq j\leq n;
$$
$$
\frac{1}{\alpha (n-j+1)^\alpha} \leq \frac{b_{jn}}{j^\alpha-(j-1)^\alpha} \leq \frac{1}{\alpha (n-j)^\alpha},\quad 1\leq j\leq n; \eqno{(3.3)}
$$
$$
\frac{1}{\alpha (n-j)^\alpha} \leq \frac{b_{j+1,n}}{(j+1)^\alpha j^\alpha} \leq \frac{1}{\alpha (n-j-1)^\alpha},\quad 1\leq j\leq n-1.
$$
By the last inequality, $a_{jn}\leq 0$, $1\leq j\leq n-1$.

Now let us estimate $|a_{0n}|$ and $a_{nn}$, $n\geq 1$. For $j=1$, the inequality (3.3) takes the form
$$
\frac{1}{\alpha n^\alpha}\leq b_{1n}\leq \frac{1}{\alpha(n-1)^\alpha}.
$$
Multiplying all parts of this inequality by $\alpha/\Gamma(1-\alpha)$, we come to the estimate for $|a_{0n}|$:
$$
\frac{1}{\Gamma(1-\alpha)n^\alpha}\leq |a_{0n}|\leq \frac{1}{\Gamma(1-\alpha)(n-1)^\alpha}. \eqno{(3.4)}
$$
Let us turn to $a_{nn}$. We have
$$
\frac{1}{n^{1-\alpha}} \int\limits_{n-1}^n \frac{dx}{(n-x)^\alpha} \leq \int\limits_{n-1}^n \frac{dx}{x^{1-\alpha}(n-x)^\alpha} \leq \frac{1}{(n-1)^{1-\alpha}} \int\limits_{n-1}^n \frac{dx}{(n-x)^\alpha}.
$$
Consequently,
$$
\frac{1}{(1-\alpha)n^{1-\alpha}}\leq b_{nn}\leq \frac{1}{(1-\alpha)(n-1)^{1-\alpha}}.
$$
In addition,
$$
\frac{\alpha}{n^{1-\alpha}}\leq n^\alpha-(n-1)^\alpha \leq \frac{\alpha}{(n-1)^{1-\alpha}}.
$$
Combining the last inequality with (3.1), we obtain the desired estimate for $a_{nn}$:
$$
\frac{1}{\Gamma(2-\alpha)}\Bigl(\frac{n-1}{n}\Bigr)^{1-\alpha} \leq a_{nn}\leq \frac{1}{\Gamma(2-\alpha)}\Bigl(\frac{n}{n-1}\Bigr)^{1-\alpha}. \eqno{(3.5)}
$$

Our next goal is to obtain estimates for $|v_n|$, $1\leq n\leq N$, in terms of $\lambda$, $\Delta t$, $n$. We observe that the current value $v_n$ can be expressed by (3.2) as
$$
v_n=-\frac{1}{a_{nn}-\lambda(\Delta t)^\alpha}\sum\limits_{j=0}^{n-1}a_{jn}v_j,\quad 1\leq n\leq N. $$
Hence,
$$
|v_n|\leq\frac{1}{|a_{nn}-\lambda(\Delta t)^\alpha|} \sum\limits_{j=0}^{n-1} |a_{jn}||v_j|, \quad 1\leq n\leq N.
$$
Since $a_{nn}>0$, by Condition 2.1, we get
$$
|a_{nn}-\lambda(\Delta t)^\alpha|\geq\sqrt{a_{nn}^2+L^2+2a_{nn}L\cos\varphi_0}, \eqno{(3.6)}
$$
where
$$
L=L(\lambda,\Delta t)=|\lambda|(\Delta t)^\alpha,
$$
and the angle $\varphi_0$ is defined in Condition 1.1. Therefore,
$$
|v_n|\leq\frac{1}{\sqrt{a_{nn}^2+L^2+2a_{nn}L\cos\varphi_0}} \sum\limits_{j=0}^{n-1} |a_{jn}||v_j|, \quad 1\leq n\leq N; \quad v_0=1. \eqno{(3.7)}
$$

It is easy to prove by induction that for any $n\geq 1$ there exists $K^*(n)>0$ such that
$$
|v_n|\leq \frac{K^*(n)}{1+L}\quad \forall L>0. \eqno{(3.8)}
$$
Indeed, we have
$$
|v_1|\leq\frac{|a_{01}|}{\sqrt{a_{11}^2+L^2+2a_{11}L\cos\varphi_0}}\leq \frac{K^*(1)}{1+L}
$$
with a constant $K^*(1)$. If $K^*(j)$, $1\leq j\leq n-1$ are already defined, then from (3.7) it follows that
$$
|v_n|\leq \frac{1}{\sqrt{a_{nn}^2+L^2+2a_{nn}L\cos\varphi_0}}\Bigl( |a_{0n}|+\sum\limits_{j=1}^{n-1}\frac{|a_{jn}|K^*(j)}{1+L} \Bigr).
$$
Hence (3.8) is fulfilled with some constant $K^*(n)$.

We now turn to more accurate estimates of $|v_n|$. For this purpose we need auxiliary inequalities (4.1) and (4.14), see Section 4. These inequalities are true for $n\geq n_0$ where $n_0=n_0(\alpha,\varepsilon)\geq 2$ and $\varepsilon>0$ is sufficiently small.

Firstly, we consider the case $\alpha\in[1/2,1)$. We claim that for sufficiently small $\varepsilon>0$,
$$
|v_n|\leq \frac{K}{Ln^{1-\alpha-\varepsilon}}, \quad n\geq 1,\quad K=K(\alpha,\varepsilon), \eqno{(3.9)}
$$
uniformly in $L>0$. According to (3.8), inequality (3.9) holds for $1\leq n\leq n_0-1$ and $K\geq \max\limits_{1\leq j\leq n_0-1}K^*(j)j^{1-\alpha-\varepsilon}$. We will prove that there exists a constant $K\geq \max\limits_{1\leq j\leq n_0-1}K^*(j)j^{1-\alpha-\varepsilon}$, such that the inductive assumption
$$
|v_j|\leq\frac{K}{Lj^{1-\alpha-\varepsilon}},\quad 1\leq j\leq n-1 \eqno{(3.10)}
$$
with any particular $n\geq n_0$ yields (3.9) with this $n$. Then the inequality (3.9) will be justified for all $n\geq 1$. From (3.7) and (3.10) it follows that
$$
|v_n|\leq \frac{1}{\sqrt{a_{nn}^2+L^2+2a_{nn}L\cos\varphi_0}} \Bigl( |a_{0n}|+\sum\limits_{j=1}^{n-1}\frac{|a_{jn}|}{j^{1-\alpha-\varepsilon}}\frac{K}{L} \Bigr).
$$
Since $n\geq n_0$, we can use the estimate (4.1) from Lemma 4.1:
$$
|v_n|\leq \frac{1}{\sqrt{1+(2L\cos\varphi_0)/a_{nn}+L^2/a_{nn}^2}} \Bigl( \frac{|a_{0n}|}{a_{nn}}+ \frac{K}{Ln^{1-\alpha-\varepsilon}} \Bigr).
$$
Taking into account (3.4), (3.5), we obtain that (3.9) holds with $K\geq \max\limits_{1\leq j\leq n_0-1}K^*(j)j^{1-\alpha-\varepsilon}$ such that for all $n\geq n_0$ and any $L>0$,
$$
\frac{1}{\sqrt{1+(2L\cos\varphi_0)/a_{nn}+L^2/a_{nn}^2}} \Bigl( (1-\alpha)\frac{n^{2(1-\alpha)-\varepsilon}}{n-1}+\frac{K}{L} \Bigr) \leq \frac{K}{L}.
$$
We can rewrite this condition on $K$ as
$$
K\geq\frac{(1-\alpha)L}{\sqrt{1+(2L\cos\varphi_0)/a_{nn}+L^2/a_{nn}^2}-1}\frac{n^{2(1-\alpha)-\varepsilon}}{n-1},\quad L>0,\quad n\geq n_0.
$$
Since $\alpha\geq 1/2$, by (3.5), we get the inequality
$$
\frac{n^{2(1-\alpha)-\varepsilon}}{n-1}\leq\frac{n_0^{2(1-\alpha)-\varepsilon}}{n_0-1},\quad a_{nn}\leq \frac{1}{\Gamma(2-\alpha)}\Bigl(\frac{n_0}{n_0-1}\Bigr)^{1-\alpha},\quad n\geq n_0.
$$
Therefore, the inductive assumption yields the estimate (3.9) with any constant $K\geq \max\limits_{1\leq j\leq n_0-1}K^*(j)j^{1-\alpha-\varepsilon}$ satisfying the condition
$$
K\geq \frac{(1-\alpha)L}{\sqrt{1+\Bigl(2\Gamma(2-\alpha)\Bigl(\frac{n_0-1}{n_0}\Bigr)^{1-\alpha}\cos\varphi_0\Bigr)L +\Gamma^2(2-\alpha)\Bigl(\frac{n_0-1}{n_0}\Bigr)^{2-2\alpha}L^2}-1} \frac{n_0^{2(1-\alpha)-\varepsilon}}{n_0-1}\,\forall L>0.
$$
The existence of such constants follows by the boundedness of the right part of the last inequality on the ray $L\in(0,+\infty)$. Thus, the estimate (3.9) is proved for any $n\geq 1$. Letting $L=L(\lambda,\Delta t)=|\lambda|(\Delta t)^\alpha$, we come to the following assertion.

{\bf Lemma 3.1.} {\it Let Condition 2.1 be fulfilled. If $\alpha\in[1/2,1)$ and $\varepsilon>0$, then for the solution of the difference equation (3.2), the estimate is true:
$$
|v_n|\leq\frac{C_{20}}{|\lambda|(\Delta t)^{\alpha}n^{1-\alpha-\varepsilon}},\quad 1\leq n\leq N.
$$
}

\medskip

Let us consider the case $\alpha\in(0,1/2)$. Using (3.7), we prove by induction that
$$
|v_n|\leq\frac{K}{Ln^\alpha},\quad n\geq 1;\quad K=K(\alpha). \eqno{(3.11)}
$$
As in the case of $\alpha\in[1/2,1)$, to establish (3.11) for some $K\geq \max\limits_{1\leq j\leq n_0-1}K^*(j)j^\alpha$, it is sufficient to prove that the inductive assumption
$$
|v_j|\leq\frac{K}{Lj^\alpha},\quad 1\leq j\leq n-1 \eqno{(3.12)}
$$
with any $n\geq n_0$ yields (3.11) with this $n$. For any $n\geq n_0$, from (3.7) and (3.12) we obtain
$$
|v_n|\leq \frac{1}{\sqrt{a_{nn}^2+L^2+2a_{nn}L\cos\varphi_0}} \Bigl( |a_{0n}|+\sum\limits_{j=1}^{n-1}\frac{|a_{jn}|}{j^\alpha}\frac{K}{L} \Bigr).
$$
It suffices to prove that
$$
\frac{1}{\sqrt{a_{nn}^2+L^2+2a_{nn}L\cos\varphi_0}} \Bigl( |a_{0n}|+\sum\limits_{j=1}^{n-1}\frac{|a_{jn}|}{j^\alpha}\frac{K}{L} \Bigr) \leq \frac{K}{Ln^\alpha},\quad n\geq n_0.
$$
This inequality holds if
$$
\frac{1}{\sqrt{1+(2L\cos\varphi_0)/a_{nn}+L^2/a_{nn}^2}} \Bigl( (1-\alpha)\frac{n}{n-1}+\frac{K}{L} \Bigr)\leq \frac{K}{L},\quad n\geq n_0.
$$
Here, we have used Corollary 4.1 and estimates (3.4), (3.5). The last inequality is satisfied with any constant $K\geq \max\limits_{1\leq j\leq n_0-1}K^*(j)j^\alpha$ for which
$$
K\geq \frac{(1-\alpha)L}{\sqrt{1+\Bigl(2\Gamma(2-\alpha)\Bigl(\frac{n_0-1}{n_0}\Bigr)^{1-\alpha}\cos\varphi_0\Bigr)L +\Gamma^2(2-\alpha)\Bigl(\frac{n_0-1}{n_0}\Bigr)^{2-2\alpha}L^2}-1} \frac{n_0}{n_0-1}\,\forall L>0.
$$
Obviously, such a constant exists. We come to the following assertion.

{\bf Lemma 3.2.} {\it Let Condition 2.1 be fulfilled. If $\alpha\in(0,1/2)$, then for the solution of the difference equation (3.2), the estimate holds:
$$
|v_n|\leq\frac{C_{21}}{|\lambda|(\Delta t)^{\alpha}n^\alpha},\quad 1\leq n\leq N.
$$
}

\medskip

With Lemmas 3.1 and 3.2, we have fulfilled one of the main tasks of this section. Now we turn to the rate--of--convergence estimates of the scheme (3.2). We denote
$$
g_n=\sum\limits_{j=0}^n \frac{a_{jn}v(t_j)}{(\Delta t)^\alpha}-\lambda v(t_n),\quad 1\leq n\leq N. \eqno{(3.13)}
$$
By Lemma 2.3,
\begin{equation}
\begin{split}
&|g_n|=\biggl| \sum\limits_{j=0}^n \frac{a_{jn}v(t_j)}{(\Delta t)^\alpha}-\lambda v(t_n) \biggr|= \biggl| \sum\limits_{j=0}^n \frac{a_{jn}v(t_j)}{(\Delta t)^\alpha}-\partial^\alpha v(t_n) \biggr| \leq\\
&\leq |\Delta_1(\lambda,n,\Delta t)+\Delta_2(\lambda,n,\Delta t)|\leq C_{19}\eta(n)|\lambda|^2(\Delta t)^\alpha(1+|\lambda|(\Delta t)^\alpha).
\end{split}
\tag{3.14}
\end{equation}
Now we subtract (3.13) from (3.2):
$$
\sum\limits_{j=0}^n \frac{a_{jn}(v_j-v(t_j))}{(\Delta t)^\alpha}-\lambda(v_n-v(t_n))=-g_n,\quad 1\leq n\leq N. \eqno{(3.15)}
$$
Below we use the notation $R_n=v_n-v(t_n)$, $0\leq n\leq N$; in particular, $R_0=0$. By (3.15) we get the recurrent equation for $R_n$:
$$
\sum\limits_{j=0}^n \frac{a_{jn}R_j}{(\Delta t)^\alpha}-\lambda R_n=-g_n,\quad 1\leq n\leq N;\quad R_0=0.
$$
It follows that
$$
R_n=\frac{1}{a_{nn}-\lambda(\Delta t)^\alpha}\biggl(-\sum\limits_{j=0}^{n-1}a_{jn}R_j - g_n(\Delta t)^\alpha\biggr),\quad 1\leq n\leq N.
$$

Let us estimate $|R_n|$ as a function of $\lambda$, $\Delta t$, $n$. Obviously,
$$
|R_n|\leq \frac{1}{|a_{nn}-\lambda(\Delta t)^\alpha|} \biggl( \sum\limits_{j=0}^{n-1}|a_{jn}||R_j| + |g_n|(\Delta t)^\alpha \biggr),\quad 1\leq n\leq N.
$$
Using (3.14), (2.18) and (3.6), we obtain
$$
|R_n|\leq \frac{1}{\sqrt{a_{nn}^2+L^2+2a_{nn}L\cos\varphi_0}} \biggl( \sum\limits_{j=0}^{n-1}|a_{jn}||R_j| + \frac{M}{n^{s(\alpha)}}\biggr),\quad 1\leq n\leq N. \eqno{(3.16)}
$$
Here,
$$
M=M(\lambda,\Delta t)=C_{22}|\lambda|^2(\Delta t)^{2\alpha}(1+|\lambda|(\Delta t)^\alpha),
$$
\begin{equation}
\begin{split}
s(\alpha)=
\begin{Bmatrix}
\alpha,&\alpha\in(0,1/2)\\
1-\alpha-\varepsilon,&\alpha\in[1/2,1)
\end{Bmatrix},
\end{split}
\tag{3.17}
\end{equation}
and $\varepsilon>0$ can be taken arbitrarily small.

We will prove the estimate
$$
|R_n|\leq\frac{\mathcal{K}M}{Ln^{s(\alpha)}},\quad n\geq 1 \eqno{(3.18)}
$$
with a constant $\mathcal{K}=\mathcal{K}(\alpha)$, independent of $M$, $L>0$. As in the proofs of Lemmas 3.1, 3.2, we apply induction by $n\geq n_0$. Assume that
$$
|R_j|\leq\frac{\mathcal{K}M}{Lj^{s(\alpha)}},\quad 1\leq j\leq n-1,\quad n\geq n_0 \eqno{(3.19)}
$$
By Lemma 4.1 and Corollary 4.1, for any $\alpha\in(0,1)$,
$$
\sum\limits_{j=1}^{n-1} \frac{|a_{jn}|}{j^{s(\alpha)}}\leq \frac{a_{nn}}{n^{s(\alpha)}},\quad n\geq n_0.
$$
Substituting this estimate together with (3.19) into (3.16), we obtain
$$
|R_n|\leq \frac{M}{\sqrt{a_{nn}^2+L^2+2a_{nn}L\cos\varphi_0}} \biggl( \frac{\mathcal{K}}{L}\frac{a_{nn}}{n^{s(\alpha)}} + \frac{1}{n^{s(\alpha)}}\biggr),\quad 1\leq n\leq N.
$$
Now to prove (3.18) it suffices to establish the existence of a constant $\mathcal{K}$ such that for any $M$, $L>0$ and $n\geq 1$,
$$
\frac{1}{\sqrt{a_{nn}^2+L^2+2a_{nn}L\cos\varphi_0}} \biggl( a_{nn}\frac{\mathcal{K}}{L} + 1\biggr) \leq \frac{\mathcal{K}}{L}
$$
or, equivalently,
$$
\mathcal{K}\geq \frac{L}{a_{nn}(\sqrt{1+(2L\cos\varphi_0)/a_{nn}+L^2/a_{nn}^2}-1)}.
$$
The existence of such constants $\mathcal{K}$ is ensured by the estimate (3.5) and the boundedness of the right side in the last inequality for $L\in(0,+\infty)$. Thus, the estimate (3.18) is proved. We get the following result.

{\bf Lemma 3.3.} {\it Let Condition 2.1 be fulfilled. Then for the difference scheme (3.2), the estimate holds:
\begin{equation}
|v_n-v(t_n)|\leq\frac{C_{23}|\lambda|(\Delta t)^\alpha(1+|\lambda|(\Delta t)^\alpha)}{n^{s(\alpha)}}
,\quad 1\leq n\leq N.
\notag
\end{equation}
Here, $s(\alpha)$ is defined in (3.17), and $C_{23}=C_{23}(\alpha,\varepsilon)$.
}

\medskip

Estimates of Lemmas 3.1--3.3 are the main result of this section. In Section 5, we will use these estimates to justify a finite difference method for solving (1.1).

\begin{center}
\textbf{4. Auxiliary inequalities}
\end{center}

This section contains the proofs of auxiliary inequalities, see Lemma 4.1 and Corollary 4.1. We have used these results in Section 3.

{\bf Lemma 4.1.} {\it For any $\alpha\in (0,1)$ and sufficiently small $\varepsilon>0$, there exists a number $n_0=n_0(\alpha,\varepsilon)$ such that
$$
\sum\limits_{j=1}^{n-1}\frac{|a_{jn}|}{j^{1-\alpha-\varepsilon}}\leq\frac{a_{nn}}{n^{1-\alpha-\varepsilon}},\quad n\geq n_0. \eqno{(4.1)}
$$}

{\bf Proof.} We rewrite the inequality (4.1) as follows:
$$
\frac{\alpha}{\Gamma(1-\alpha)} \sum\limits_{j=1}^{n-1} \Bigl(\frac{b_{j+1,n}}{j^{1-\alpha-\varepsilon}\bigl((j+1)^\alpha j^\alpha\bigr)} -\frac{b_{jn}}{j^{1-\alpha-\varepsilon}\bigl(j^\alpha-(j-1)^\alpha\bigr)}\Bigr) \leq
$$
$$
\leq\frac{\alpha}{\Gamma(1-\alpha)}\frac{b_{nn}}{n^{1-\alpha-\varepsilon}\bigl(n^\alpha-(n-1)^\alpha\bigr)},\quad n\geq n_0;
$$
$$
\sum\limits_{j=2}^n \Bigl( \frac{1}{(j-1)^{1-\alpha-\varepsilon}}-\frac{1}{j^{1-\alpha-\varepsilon}} \Bigr)\frac{b_{jn}}{j^\alpha-(j-1)^\alpha} \leq b_{1n}.
$$
Here, we have used formula (3.1) and the fact that $a_{jn}\leq 0$, $0\leq j\leq n-1$. Let us now present the desired inequality as
$$
\sum\limits_{j=2}^n \Bigl( \frac{1}{(j-1)^{1-\alpha-\varepsilon}}-\frac{1}{j^{1-\alpha-\varepsilon}} \Bigr)\frac{1}{j^\alpha-(j-1)^\alpha} \int\limits_{\frac{j-1}{n}}^{\frac{j}{n}} \frac{dx}{x^{1-\alpha}(1-x)^\alpha} \leq \int\limits_{0}^{\frac{1}{n}} \frac{dx}{x^{1-\alpha}(1-x)^\alpha}.
\eqno{(4.2)}
$$

Below we need the incomplete Beta-function
$$
B_x(p,q)=\int\limits_0^x t^{p-1}(1-t)^{q-1}dt
$$
and hypergeometric functions
$$
F(a,b;c;z)={}_2F_1(a,b;c;z)=\sum\limits_{n=0}^{\infty}\frac{(a)_n(b)_n}{(c)_n n!}z^n, \quad (a)_n=\frac{\Gamma(a+n)}{\Gamma(a)},\quad a,b,c\neq 0,-1,-2,\dots .
$$
The series in the definition of $F(a,b;c;z)$ converges absolutely when $|z|<1$, and when $|z|=1$ if ${\rm{Re}}(a+b-c)<0$. We recall the identity [20,Ch.2]
$$
B_x(p,q)=p^{-1}x^p F(p,1-q;p+1;x). \eqno{(4.3)}
$$

We are now ready to prove (4.2). Using (4.3) we get
$$
\int\limits_{\frac{j-1}{n}}^{\frac{j}{n}} \frac{dx}{x^{1-\alpha}(1-x)^\alpha} = B_{\frac{j}{n}}(\alpha,1-\alpha)-B_{\frac{j-1}{n}}(\alpha,1-\alpha) =
$$
$$
=\frac{1}{\alpha n^\alpha} \Bigl(j^\alpha \mathcal{F}\Bigl(\frac{j}{n}\Bigr)-(j-1)^\alpha \mathcal{F}\Bigl(\frac{j-1}{n}\Bigr)\Bigr),\quad 1\leq j\leq n,
$$
where $\mathcal{F}(x)=F(\alpha,\alpha;1+\alpha;x)$. We rewrite (4.2) as follows:
$$
\sum\limits_{j=2}^n \Bigl( \frac{1}{(j-1)^{1-\alpha-\varepsilon}}-\frac{1}{j^{1-\alpha-\varepsilon}} \Bigr)\frac{1}{j^\alpha-(j-1)^\alpha} \Bigl( j^\alpha \mathcal{F}\Bigl(\frac{j}{n}\Bigr) -(j-1)^\alpha \mathcal{F}\Bigl(\frac{j-1}{n}\Bigr) \Bigr) \leq \mathcal{F}\Bigl(\frac{1}{n}\Bigr). \eqno{(4.4)}
$$
The function $\mathcal{F}(x)$ has the form
$$
\mathcal{F}(x)=\sum\limits_{n=0}^\infty \frac{\alpha\Gamma(n+\alpha)}{(n+\alpha)n!\Gamma(\alpha)}x^n =1+x\mathcal{G}(x),\quad \mathcal{G}(x)=\frac{\alpha}{\Gamma(\alpha)} \sum\limits_{n=0}^\infty \frac{\Gamma(n+1+\alpha)}{(n+1+\alpha)(n+1)!}x^n.
$$
It is easy to see that $\mathcal{G}(x)$ is increasing when $x\geq 0$. We rewrite the left part of (4.4) as
\begin{equation}
\begin{split}
\sum\limits_{j=2}^n \Bigl( \frac{1}{(j-1)^{1-\alpha-\varepsilon}}&-\frac{1}{j^{1-\alpha-\varepsilon}} \Bigr)\frac{1}{j^\alpha-(j-1)^\alpha} \Bigl( j^\alpha \mathcal{F}\Bigl(\frac{j}{n}\Bigr) -(j-1)^\alpha \mathcal{F}\Bigl(\frac{j-1}{n}\Bigr) \Bigr) =\\
=\sum\limits_{j=2}^n \Bigl( \frac{1}{(j-1)^{1-\alpha-\varepsilon}}&-\frac{1}{j^{1-\alpha-\varepsilon}} \Bigr) +\frac{1}{n}\sum\limits_{j=2}^n \Bigl( \frac{1}{(j-1)^{1-\alpha-\varepsilon}}-\frac{1}{j^{1-\alpha-\varepsilon}} \Bigr)\frac{1}{j^\alpha-(j-1)^\alpha} \cdot\\
&\cdot\Bigl( j^{1+\alpha}\mathcal{G}\Bigl(\frac{j}{n}\Bigr) -(j-1)^{1+\alpha}\mathcal{G}\Bigl(\frac{j-1}{n}\Bigr) \Bigr)=\\
=\Bigl(1-\frac{1}{n^{1-\alpha-\varepsilon}}\Bigr) &+\frac{1}{n}\sum\limits_{j=2}^n \Bigl( \frac{1}{(j-1)^{1-\alpha-\varepsilon}}-\frac{1}{j^{1-\alpha-\varepsilon}} \Bigr)\frac{j^{1+\alpha} (j-1)^{1+\alpha}}{j^\alpha-(j-1)^\alpha} \mathcal{G}\Bigl(\frac{j}{n}\Bigr)+\\
+\frac{1}{n}\sum\limits_{j=2}^n \Bigl( \frac{1}{(j-1)^{1-\alpha-\varepsilon}}&-\frac{1}{j^{1-\alpha-\varepsilon}} \Bigr)\frac{(j-1)^{1+\alpha}}{j^\alpha-(j-1)^\alpha} \Bigl( \mathcal{G}\Bigl(\frac{j}{n}\Bigr) -\mathcal{G}\Bigl(\frac{j-1}{n}\Bigr) \Bigr).
\end{split}
\tag{4.5}
\end{equation}
Let us estimate the first sum in the right part of the last equality:
$$
\frac{1}{n}\sum\limits_{j=2}^n \Bigl( \frac{1}{(j-1)^{1-\alpha-\varepsilon}}-\frac{1}{j^{1-\alpha-\varepsilon}} \Bigr)\frac{j^{1+\alpha} (j-1)^{1+\alpha}}{j^\alpha-(j-1)^\alpha} \mathcal{G}\Bigl(\frac{j}{n}\Bigr) \leq
$$
$$
\leq\frac{1}{n}\sum\limits_{j=2}^n \bigl( (1-\alpha-\varepsilon)(j-1)^{-2+\alpha+\varepsilon} \bigr)\frac{(1+\alpha)j^\alpha}{\alpha j^{-1+\alpha}}\mathcal{G}\Bigl(\frac{j}{n}\Bigr) =\frac{(1+\alpha)(1-\alpha-\varepsilon)}{\alpha n} \sum\limits_{j=2}^n \frac{j}{(j-1)^{2-\alpha-\varepsilon}}\mathcal{G}\Bigl(\frac{j}{n}\Bigr)=
$$
$$
=\frac{(1+\alpha)(1-\alpha-\varepsilon)}{\alpha n} \sum\limits_{j=2}^n \frac{1}{j^{1-\alpha-\varepsilon}}\mathcal{G}\Bigl(\frac{j}{n}\Bigr) +\frac{(1+\alpha)(1-\alpha-\varepsilon)}{\alpha n} \sum\limits_{j=2}^n \Bigl( \frac{j}{(j-1)^{2-\alpha-\varepsilon}} -\frac{1}{j^{1-\alpha-\varepsilon}} \Bigr)\mathcal{G}\Bigl(\frac{j}{n}\Bigr).
$$
Here,
$$
\frac{(1+\alpha)(1-\alpha-\varepsilon)}{\alpha n} \sum\limits_{j=2}^n \frac{1}{j^{1-\alpha-\varepsilon}}\mathcal{G}\Bigl(\frac{j}{n}\Bigr) = \frac{(1+\alpha)(1-\alpha-\varepsilon)}{\alpha n^{1-\alpha-\varepsilon}} \cdot \frac{1}{n} \sum\limits_{j=2}^n \frac{1}{\Bigl(\frac{j}{n}\Bigr)^{1-\alpha-\varepsilon}}\mathcal{G}\Bigl(\frac{j}{n}\Bigr) \leq
$$
$$
\leq (1+o(1))\frac{(1+\alpha)(1-\alpha-\varepsilon)}{\alpha}\cdot\frac{1}{n^{1-\alpha-\varepsilon}}\int\limits_0^1 x^{-1+\alpha+\varepsilon}\mathcal{G}(x)dx;
$$
$$
\frac{(1+\alpha)(1-\alpha-\varepsilon)}{\alpha n} \sum\limits_{j=2}^n \Bigl( \frac{j}{(j-1)^{2-\alpha-\varepsilon}} -\frac{1}{j^{1-\alpha-\varepsilon}} \Bigr)\mathcal{G}\Bigl(\frac{j}{n}\Bigr) \leq
$$
$$
\leq\frac{(1+\alpha)(1-\alpha-\varepsilon)\mathcal{G}(1)}{\alpha n} \biggl( \sum\limits_{j=2}^n\Bigl( \frac{1}{(j-1)^{1-\alpha-\varepsilon}} -\frac{1}{j^{1-\alpha-\varepsilon}} \Bigr) +\sum\limits_{j=2}^n\frac{1}{(j-1)^{2-\alpha-\varepsilon}} \biggr)\leq
$$
$$
\leq \frac{(1+\alpha)(1-\alpha-\varepsilon)\mathcal{G}(1)}{\alpha n} \biggl( \Bigl(1-\frac{1}{n^{1-\alpha-\varepsilon}}\Bigr) +\biggl( 1+\int\limits_2^n\frac{dx}{(x-1)^{2-\alpha-\varepsilon}} \biggr) \biggr)=
$$
$$
=\frac{(1+\alpha)(1-\alpha-\varepsilon)\mathcal{G}(1)}{\alpha n} \Bigl( 2-\frac{1}{n^{1-\alpha-\varepsilon}}+\frac{1-(n-1)^{-1+\alpha+\varepsilon}}{1-\alpha-\varepsilon} \Bigr)\leq \frac{C_{24}}{n}.
$$
Consequently,
\begin{equation}
\begin{split}
&\frac{1}{n}\sum\limits_{j=2}^n \Bigl( \frac{1}{(j-1)^{1-\alpha-\varepsilon}}-\frac{1}{j^{1-\alpha-\varepsilon}} \Bigr)\frac{j^{1+\alpha} (j-1)^{1+\alpha}}{j^\alpha-(j-1)^\alpha} \mathcal{G}\Bigl(\frac{j}{n}\Bigr) \leq\\
&\leq (1+o(1))\frac{(1+\alpha)(1-\alpha-\varepsilon)}{\alpha}\cdot\frac{1}{n^{1-\alpha-\varepsilon}}\int\limits_0^1 x^{-1+\alpha+\varepsilon}\mathcal{G}(x)dx.
\end{split}
\tag{4.6}
\end{equation}
Here and below, $o(1)$ means an infinitesimal, as $n\to\infty$. We have
$$
\frac{1}{n}\sum\limits_{j=2}^n \Bigl( \frac{1}{(j-1)^{1-\alpha-\varepsilon}}-\frac{1}{j^{1-\alpha-\varepsilon}} \Bigr)\frac{(j-1)^{1+\alpha}}{j^\alpha-(j-1)^\alpha} \Bigl( \mathcal{G}\Bigl(\frac{j}{n}\Bigr) -\mathcal{G}\Bigl(\frac{j-1}{n}\Bigr) \Bigr)=
$$
$$
=\frac{1}{n}\sum\limits_{j=2}^n \frac{j^{1-\alpha-\varepsilon} (j-1)^{1-\alpha-\varepsilon}}{j^\alpha-(j-1)^\alpha} \cdot\frac{(j-1)^{2\alpha+\varepsilon}}{j^{1-\alpha-\varepsilon}} \Bigl( \mathcal{G}\Bigl(\frac{j}{n}\Bigr) -\mathcal{G}\Bigl(\frac{j-1}{n}\Bigr) \Bigr)\leq
$$
$$
\leq \frac{1}{n}\sum\limits_{j=2}^n \frac{(1-\alpha-\varepsilon)j^{1-\alpha}}{\alpha (j-1)^{\alpha+\varepsilon}} \cdot\frac{(j-1)^{2\alpha+\varepsilon}}{j^{1-\alpha-\varepsilon}} \Bigl( \mathcal{G}\Bigl(\frac{j}{n}\Bigr) -\mathcal{G}\Bigl(\frac{j-1}{n}\Bigr) \Bigr)=
$$
$$
=\frac{1-\alpha-\varepsilon}{\alpha n} \sum\limits_{j=2}^n (j-1)^\alpha j^\varepsilon \Bigl( \mathcal{G}\Bigl(\frac{j}{n}\Bigr) -\mathcal{G}\Bigl(\frac{j-1}{n}\Bigr) \Bigr) \leq
$$
$$
\leq\frac{2^\varepsilon(1-\alpha-\varepsilon)}{\alpha n} \sum\limits_{j=1}^n (j-1)^{\alpha+\varepsilon} \Bigl( \mathcal{G}\Bigl(\frac{j}{n}\Bigr) -\mathcal{G}\Bigl(\frac{j-1}{n}\Bigr) \Bigr)\leq \frac{2^\varepsilon(1-\alpha-\varepsilon)}{\alpha n^{1-\alpha-\varepsilon}} \int\limits_0^1 x^{\alpha+\varepsilon} d\mathcal{G}(x).
$$
It follows that
\begin{equation}
\begin{split}
&\frac{1}{n}\sum\limits_{j=2}^n \Bigl( \frac{1}{(j-1)^{1-\alpha-\varepsilon}}-\frac{1}{j^{1-\alpha-\varepsilon}} \Bigr)\frac{(j-1)^{1+\alpha}}{j^\alpha-(j-1)^\alpha} \Bigl( \mathcal{G}\Bigl(\frac{j}{n}\Bigr) -\mathcal{G}\Bigl(\frac{j-1}{n}\Bigr) \Bigr)\leq\\
&\leq \frac{2^\varepsilon(1-\alpha-\varepsilon)}{\alpha n^{1-\alpha-\varepsilon}}\biggl( \mathcal{G}(1)-(\alpha+\varepsilon) \int\limits_0^1 x^{-1+\alpha+\varepsilon}\mathcal{G}(x)dx \biggr).
\end{split}
\tag{4.7}
\end{equation}
By substituting (4.6) and (4.7) into (4.5), we get
$$
\sum\limits_{j=2}^n \Bigl( \frac{1}{(j-1)^{1-\alpha-\varepsilon}}-\frac{1}{j^{1-\alpha-\varepsilon}} \Bigr)\frac{1}{j^\alpha-(j-1)^\alpha} \Bigl( j^\alpha \mathcal{F}\Bigl(\frac{j}{n}\Bigr) -(j-1)^\alpha \mathcal{F}\Bigl(\frac{j-1}{n}\Bigr) \Bigr) \leq
$$
$$
\leq\Bigl(1-\frac{1}{n^{1-\alpha-\varepsilon}}\Bigr) +(1+o(1))\frac{(1+\alpha)(1-\alpha-\varepsilon)}{\alpha}\cdot\frac{1}{n^{1-\alpha-\varepsilon}}\int\limits_0^ 1 x^{-1+\alpha+\varepsilon}\mathcal{G}(x)dx+
$$
$$
+\frac{2^\varepsilon(1-\alpha-\varepsilon)}{\alpha}\biggl( \frac{\mathcal{G}(1)}{n^{1-\alpha-\varepsilon}} -\frac{\alpha+\varepsilon}{n^{1-\alpha-\varepsilon}}\int\limits_0^1 x^{-1+\alpha+\varepsilon}\mathcal{G}(x)dx \biggr).
$$
To determine the value of $\mathcal{G}(1)$ we observe that according to (4.3),
$$
B_1(\alpha,1-\alpha)=\frac{1}{\alpha}F(\alpha,\alpha;1+\alpha;1) =\frac{1}{\alpha}\mathcal{F}(1).
$$
On the other hand,
$$
B_1(\alpha,1-\alpha)=B(\alpha,1-\alpha)=\frac{\pi}{\sin\alpha\pi}.
$$
Thus,
$$
\mathcal{F}(1)=\frac{\alpha\pi}{\sin\alpha\pi};\quad \mathcal{G}(1)=\mathcal{F}(1)-1=\frac{\alpha\pi}{\sin\alpha\pi}-1.
$$
Therefore,
$$
\sum\limits_{j=2}^n \Bigl( \frac{1}{(j-1)^{1-\alpha-\varepsilon}}-\frac{1}{j^{1-\alpha-\varepsilon}} \Bigr)\frac{1}{j^\alpha-(j-1)^\alpha} \Bigl( j^\alpha \mathcal{F}\Bigl(\frac{j}{n}\Bigr) -(j-1)^\alpha \mathcal{F}\Bigl(\frac{j-1}{n}\Bigr) \Bigr) \leq
$$
$$
\leq 1-\Bigl( 1-\frac{2^\varepsilon(1-\alpha-\varepsilon)}{\alpha}\Bigl( \frac{\alpha\pi}{\sin\alpha\pi}-1 \Bigr) \Bigr)\frac{1}{n^{1-\alpha-\varepsilon}}+
$$
$$
+(1+o(1))\frac{1-\alpha-\varepsilon}{\alpha} \bigl( 1-\alpha(2^\varepsilon-1)-\varepsilon 2^\varepsilon \bigr) \cdot\frac{1}{n^{1-\alpha-\varepsilon}} \int\limits_0^1 x^{-1+\alpha+\varepsilon}\mathcal{G}(x)dx.
$$
Using this inequality and the relation $\mathcal{F}\Bigl(\frac{1}{n}\Bigr)=1+\frac{1}{n}\mathcal{G}\Bigl(\frac{1}{n}\Bigr)\geq 1$, we conclude that inequality (4.4) is true for all $n$, starting from some $n_0$, if
$$
1-\frac{2^\varepsilon(1-\alpha-\varepsilon)}{\alpha}\Bigl( \frac{\alpha\pi}{\sin\alpha\pi}-1 \Bigr) -\frac{1-\alpha-\varepsilon}{\alpha} \bigl( 1-\alpha(2^\varepsilon-1)-\varepsilon 2^\varepsilon \bigr) \int\limits_0^1 x^{-1+\alpha+\varepsilon}\mathcal{G}(x)dx >0. \eqno{(4.8)}
$$

Let us prove (4.8) for sufficiently small $\varepsilon>0$. We have:
\begin{equation}
\begin{split}
\int\limits_0^1 x^{-1+\alpha+\varepsilon}\mathcal{G}(x)dx =&\int\limits_0^1 x^{-2+\alpha+\varepsilon}(\mathcal{F}(x)-1)dx=\\
=\frac{1}{1-\alpha-\varepsilon}\biggl( &\int\limits_0^1 x^{-1+\alpha+\varepsilon}\mathcal{F}^\prime(x)dx -\Bigl(\frac{\alpha\pi}{\sin\alpha\pi}-1\Bigr) \biggr).
\end{split}
\tag{4.9}
\end{equation}
Applying the formula for the derivative of the hypergeometric function [20,Ch.2]
$$
\frac{d^n}{dx^n}F(a,b;c;x)=\frac{(a)_n(b)_n}{(c)_n}F(a+n,b+n;c+n;x),
$$
we obtain the following expression for $\mathcal{F}^\prime(x)$:
$$
\mathcal{F}^\prime(x)=\frac{d}{dx}F(\alpha,\alpha;1+\alpha;x)=\frac{\alpha^2}{1+\alpha}F(1+\alpha,1+\alpha;2+\alpha;x).
$$
To calculate the integral in the right part of (4.9), we use the formula [21,2.21.1(5)]:
$$
\int\limits_0^y x^{A-1}(y-x)^{B-1}F(a,b;c;\frac{x}{y})dx =B(A,B)y^{A+B-1}{}_3F_2(a,b,A;c,A+B;1).
$$
Here,
$$
{}_pF_q(a_1,\dots,a_p;b_1,\dots,b_q;x) =\sum\limits_{n=0}^\infty \frac{(a_1)_n\cdot\dots\cdot(a_p)_n}{(b_1)_n\cdot\dots\cdot(b_q)_n} \cdot\frac{x^n}{n!},
$$
and ${}_pF_q(a_1,\dots,a_p;b_1,\dots,b_q;x)$ does not depend on the order of arguments $a_1,\dots,a_p$ and $b_1,\dots,b_q$ in each of these two groups [21,7.2.3]. We get
$$
\int\limits_0^1 x^{-1+\alpha+\varepsilon}\mathcal{F}^\prime(x)dx =\frac{\alpha^2}{1+\alpha}\int\limits_0^1 x^{-1+\alpha+\varepsilon}F(1+\alpha,1+\alpha;2+\alpha;x)dx=
$$
$$
=\frac{\alpha^2}{1+\alpha}B(\alpha+\varepsilon,1){}_3F_2(1+\alpha,1+\alpha,\alpha+\varepsilon;2+\alpha,1+\alpha+\varepsilon;1)=
$$
$$
=\frac{\alpha^2}{(\alpha+\varepsilon)(1+\alpha)}{}_3F_2(1+\alpha,\alpha+\varepsilon,1+\alpha;2+\alpha,1+\alpha+\varepsilon;1).
$$
Now let us apply formula [21,7.4.4(16)]
$$
{}_3F_2(a,b,c;a+1,b+1;1)= \frac{ab}{a-b}\Gamma(1-c)\Bigl( \frac{\Gamma(b)}{\Gamma(1+b-c)}-\frac{\Gamma(a)}{\Gamma(1+a-c)} \Bigr),\quad a\neq b,c\neq 1,{\rm{Re}}c<2.
$$
With its help, we find
$$
\int\limits_0^1 x^{-1+\alpha+\varepsilon}\mathcal{F}^\prime(x)dx =\frac{\alpha^2}{(\alpha+\varepsilon)(1+\alpha)} \cdot\frac{(1+\alpha)(\alpha+\varepsilon)}{1-\varepsilon}\Gamma(-\alpha) \Bigl( \frac{\Gamma(\alpha+\varepsilon)}{\Gamma(\varepsilon)} - \frac{\Gamma(1+\alpha)}{\Gamma(1)} \Bigr)=
$$
$$
=\frac{\alpha}{1-\varepsilon}\Gamma(1-\alpha)\Bigl( \alpha\Gamma(\alpha) -\frac{\Gamma(\alpha+\varepsilon)}{\Gamma(\varepsilon)} \Bigr) =\frac{\alpha}{1-\varepsilon}\Bigl( \frac{\alpha\pi}{\sin\alpha\pi} -\frac{\Gamma(\alpha+\varepsilon)\Gamma(1-\alpha)}{\Gamma(\varepsilon)} \Bigr).
$$
Combining this with (4.9) we obtain
$$
\int\limits_0^1 x^{-1+\alpha+\varepsilon}\mathcal{G}(x)dx =\frac{1}{1-\alpha-\varepsilon}\Bigl( 1-\frac{\alpha}{1-\varepsilon}\frac{\Gamma(\alpha+\varepsilon)\Gamma(1-\alpha)}{\Gamma(\varepsilon)} -\frac{\alpha\pi}{\sin\alpha\pi}\Bigl( 1-\frac{\alpha}{1-\varepsilon} \Bigr) \Bigr).
$$
Now the desired inequality (4.8) takes the form
\begin{equation}
\begin{split}
&1-\frac{2^\varepsilon(1-\alpha-\varepsilon)}{\alpha}\Bigl( \frac{\alpha\pi}{\sin\alpha\pi}-1 \Bigr)-\\
-\frac{1}{\alpha}\bigl( 1-\alpha(2^\varepsilon-1)-\varepsilon 2^\varepsilon \bigr) \Bigl( &1-\frac{\alpha}{1-\varepsilon}\frac{\Gamma(\alpha+\varepsilon)\Gamma(1-\alpha)}{\Gamma(\varepsilon)} -\frac{\alpha\pi}{\sin\alpha\pi}\Bigl( 1-\frac{\alpha}{1-\varepsilon} \Bigr) \Bigr)>0.
\end{split}
\tag{4.10}
\end{equation}
Simple calculations prove that $\sin\alpha\pi>\alpha(1-\alpha)\pi$, $\alpha\in(0,1)$. Further, when $\varepsilon\to +0$ we have
$$
\alpha(2^\varepsilon-1)+\varepsilon 2^\varepsilon\to +0,\quad \frac{\alpha}{1-\varepsilon} \frac{\Gamma(\alpha+\varepsilon)\Gamma(1-\alpha)}{\Gamma(\varepsilon)}\to +0,\quad \frac{\alpha\pi}{\sin\alpha\pi}\Bigl( 1-\frac{\alpha}{1-\varepsilon}\Bigr)\to \frac{\alpha(1-\alpha)\pi}{\sin\alpha\pi}<1.
$$
Taking into account these relations and the inequality $\Gamma(\varepsilon)=\Gamma(1+\varepsilon)/\varepsilon\leq 1/\varepsilon$, we conclude that to prove (4.10) it suffices to show that for small $\varepsilon>0$,
$$
1-\frac{2^\varepsilon(1-\alpha-\varepsilon)}{\alpha}\Bigl( \frac{\alpha\pi}{\sin\alpha\pi}-1 \Bigr) -\frac{1}{\alpha}\Bigl( 1-\frac{\alpha\varepsilon}{1-\varepsilon}\Gamma(\alpha+\varepsilon)\Gamma(1-\alpha) -\frac{\alpha\pi}{\sin\alpha\pi}\Bigl( 1 -\frac{\alpha}{1-\varepsilon} \Bigr) \Bigr)>0.
$$
Equivalently,
\begin{equation}
\begin{split}
&\frac{2^\varepsilon(1-\alpha-\varepsilon)-(1-\alpha)}{\alpha} +\frac{\varepsilon}{1-\varepsilon}\Gamma(\alpha+\varepsilon)\Gamma(1-\alpha)+\\ &+\frac{\alpha\pi}{\sin\alpha\pi}\Bigl( \frac{1-2^\varepsilon(1-\alpha-\varepsilon)}{\alpha} -\frac{1}{1-\varepsilon} \Bigr)>0.
\end{split}
\tag{4.11}
\end{equation}

We get
$$
\frac{2^\varepsilon(1-\alpha-\varepsilon)-(1-\alpha)}{\alpha} +\frac{\varepsilon}{1-\varepsilon}\Gamma(\alpha+\varepsilon)\Gamma(1-\alpha) \geq
$$
$$
\geq\frac{(1+\varepsilon\ln 2)(1-\alpha-\varepsilon)-(1-\alpha)}{\alpha} +\frac{\varepsilon}{1-\varepsilon}\Gamma(\alpha+\varepsilon)\Gamma(1-\alpha)=
$$
$$
=\varepsilon\Bigl( \frac{\Gamma(\alpha+\varepsilon)\Gamma(1-\alpha)}{1-\varepsilon} -\frac{1}{\alpha}\bigl( 1-(1-\alpha-\varepsilon)\ln 2 \bigr) \Bigr) \geq \varepsilon\Bigl( \frac{\Gamma(\alpha+\varepsilon)\Gamma(1-\alpha)}{1-\varepsilon} -\frac{1}{\alpha}\Bigr)>0
$$
for sufficiently small $\varepsilon>0$. Indeed, as $\varepsilon\to 0$, we have
$$
\frac{\Gamma(\alpha+\varepsilon)\Gamma(1-\alpha)}{1-\varepsilon} -\frac{1}{\alpha}\, \to\, \Gamma(\alpha)\Gamma(1-\alpha) -\frac{1}{\alpha} =\frac{1}{\alpha}\Bigl( \frac{\alpha\pi}{\sin\alpha\pi}-1 \Bigr)>0.
$$
Thus,
$$
\frac{2^\varepsilon(1-\alpha-\varepsilon)-(1-\alpha)}{\alpha} +\frac{\varepsilon}{1-\varepsilon}\Gamma(\alpha+\varepsilon)\Gamma(1-\alpha)>0. \eqno{(4.12)}
$$
For small $\varepsilon>0$ it is also true that
$$
\frac{1-2^\varepsilon(1-\alpha-\varepsilon)}{\alpha} -\frac{1}{1-\varepsilon}>0. \eqno{(4.13)}
$$
Indeed, for small $\varepsilon>0$ we have $2^\varepsilon<1+\varepsilon$, therefore
$$
(1-\varepsilon)(1-2^\varepsilon(1-\alpha-\varepsilon))> (1-\varepsilon)(1-(1+\varepsilon)(1-\varepsilon-\alpha)) =
$$
$$
=\alpha(1-\varepsilon^2)+\varepsilon^2(1-\varepsilon)> \alpha(1-\varepsilon^2)+\varepsilon^2\alpha =\alpha.
$$
The inequality (4.11) follows directly by (4.12) and (4.13). By (4.11) we immediately get (4.10), (4.8), (4.4) and (4.1). The lemma is proved.

{\bf Remark 4.1.} Computer calculations confirm that inequality (4.1) in Lemma 4.1 holds true for all $n\geq 2$ for any $\alpha\in(0,1)$ and $\varepsilon\in(0,1-\alpha)$. In other words, one can let $n_0(\alpha,\varepsilon)\equiv 2$.

{\bf Corollary 4.1.} {\it If $\alpha\in(0,1/2)$, then with a sufficiently small $\varepsilon>0$,
$$
\sum\limits_{j=1}^{n-1}\frac{|a_{jn}|}{j^\alpha}\leq\frac{a_{nn}}{n^\alpha},\quad n\geq n_0(\alpha,\varepsilon). \eqno{(4.14)}
$$}

In fact, if $\alpha\in(0,1/2)$ and $\varepsilon>0$ is sufficiently small, then according to Lemma 4.1,
$$
\sum\limits_{j=1}^{n-1}\frac{|a_{jn}|}{j^\alpha} \leq n^{1-2\alpha-\varepsilon}\sum\limits_{j=1}^{n-1}\frac{|a_{jn}|}{j^{1-\alpha-\varepsilon}} \leq n^{1-2\alpha-\varepsilon}\frac{a_{nn}}{n^{1-\alpha-\varepsilon}} \leq\frac{a_{nn}}{n^\alpha}.
$$
This completes the proof of the corollary.

\begin{center}
\textbf{5. Finite difference method for solving Cauchy problems\\ with the Caputo fractional derivative in a Banach space}
\end{center}

Let us return to the problem (1.1). Assume that the operator $A$ satisfies Condition 1.1. To approximate the solution of (1.1), we use a finite difference scheme similar to (3.2):
$$
\sum\limits_{j=0}^n \frac{a_{jn}u_j}{(\Delta t)^\alpha}-Au_n=0,\quad 1\leq n\leq N;\quad u_0=f. \eqno{(5.1)}
$$
Here $u_n\in X$ is the required approximation to $u(t_n)$, $0\leq n\leq N$. The purpose of this section is to study the convergence of the scheme (5.1) and to obtain its rate--of--convergence estimate.

Along with the scheme (5.1), we consider its scalar analogue (3.2) which gives $v_n=v_n(\lambda)$ for $\lambda\in -\overline{K(\varphi_0)}$. We define the functions $w_n(\zeta)=v_n(-\zeta)$, $\zeta\in \overline{K(\varphi_0)}$. It is easy to see that
$$
w_n(\zeta)=-\frac{1}{a_{nn}+\zeta(\Delta t)^\alpha}\sum\limits_{j=0}^{n-1} a_{jn}w_j(\zeta),\quad w_0(\zeta)\equiv 1. \eqno{(5.2)}
$$
The elements of $u_n$, $0\leq n\leq N$, satisfy the similar equality
$$
u_n=-(a_{nn}E-(\Delta t)^\alpha A)^{-1}\sum\limits_{j=0}^{n-1}a_{jn}u_j,\quad u_0=f. \eqno{(5.3)}
$$
Let us prove that $u_n=w_n(-A)f$, $n\geq 1$, in the sense of the calculus of sectorial operators. We have
$$
w_1(\zeta)=-\frac{a_{01}}{a_{11}+\zeta(\Delta t)^\alpha},\quad u_1=-a_{01}(a_{11}E-(\Delta t)^\alpha A)^{-1}f=w_1(-A)f.
$$
Suppose that $u_j=w_j(-A)f$, $1\leq j\leq n-1$. Then by (5.2) and (5.3) it follows that
$$
u_n =-(a_{nn}E-(\Delta t)^\alpha A)^{-1}\biggl( a_{0n}E+\sum\limits_{j=1}^{n-1}a_{jn}w_j(-A) \biggr)f =w_n(-A)f.
$$
Thus, we see that
$$
u_n=w_n(-A)f= \frac{1}{2\pi i} \int\limits_{\Gamma(r_0,\varphi_0)} v_n(-\zeta) R(\zeta-A)fd\zeta, \quad 1\leq n\leq N.
$$
According to Lemma 1.1,
$$
u(t)=F_t(-A)f= \frac{1}{2\pi i} \int\limits_{\Gamma(r_0,\varphi_0)} E_\alpha(-\zeta t^\alpha)R(\zeta,-A)fd\zeta.
$$
Combining the last two equalities, we get the representation
$$
u_n-u(t_n)=\frac{1}{2\pi i} \int\limits_{\Gamma(r_0,\varphi_0)} \bigl( v_n(-\zeta) -E_\alpha(-\zeta t_n^\alpha) \bigr)R(\zeta,-A)fd\zeta =
$$
$$
= \frac{1}{2\pi i} \int\limits_{\Gamma(\varphi_0)} \bigl( v_n(-\zeta) -E_\alpha(-\zeta t_n^\alpha) \bigr)R(\zeta,-A)fd\zeta,\quad 1\leq n\leq N,
$$
where $\Gamma(\varphi_0)$ is the boundary of $K(\varphi_0)$. Using the change of variable $\zeta=-\lambda$ and the identity $R(-\lambda,-A)=-R(\lambda,A)$, we obtain
$$
u_n-u(t_n)=\frac{1}{2\pi i}\int\limits_{-\Gamma(\varphi_0)} \bigl( v_n(\lambda)-E_\alpha(\lambda t_n^\alpha) \bigr) R(\lambda,A)fd\lambda,\quad 1\leq n\leq N.
$$

It is suitable to divide the contour $\Gamma(\varphi_0)=\{\lambda\in\mathbb{C}\,|\,-\lambda\in\Gamma(\varphi_0)\}$ into two parts, namely, $\Gamma_1=\{\lambda\in-\Gamma(\varphi_0)\,|\,|\lambda|\leq(\Delta t)^{-\alpha}\}$ and $\Gamma_2=\{\lambda\in-\Gamma(\varphi_0)\,|\,|\lambda|\geq(\Delta t)^{-\alpha}\}$.
We have the following estimate for the error of the scheme (5.1) at $t=t_n$, $1\leq n\leq N$:
$$
\|u_n-u(t_n)\|\leq C_{25}\|f\|\biggl( \int\limits_{\Gamma_1}\frac{|v_n(\lambda)-v(t_n)|}{1+|\lambda|}|d\lambda| +\int\limits_{\Gamma_2}\frac{|v_n(\lambda)|}{1+|\lambda|}|d\lambda| +\int\limits_{\Gamma_2} \frac{|E_\alpha(\lambda t_n^\alpha)|}{1+|\lambda|}|d\lambda| \biggr). \eqno{(5.4)}
$$
Here, we have used the estimate
$$
\|R(\lambda,A)\|=\|R(-\lambda-A)\|\leq\frac{C_0}{1+|\lambda|}, \quad \lambda\in-\Gamma(\varphi_0)
$$
and the fact that $v(t)=E_\alpha(\lambda t^\alpha)$ is the solution of (2.1).

Let us estimate the first term in the right part of (5.4). Applying Lemma 3.3 we get
$$
\int\limits_{\Gamma_1}\frac{|v_n(\lambda)-v(t_n)|}{1+|\lambda|}|d\lambda| \leq \frac{C_{23} (\Delta t) ^{\alpha}}{n^{s(\alpha)}} \int\limits_{\Gamma_1} \frac{|\lambda|(1+|\lambda|(\Delta t)^\alpha)}{1+|\lambda|} |d\lambda|\leq
$$
$$
\leq \frac{C_{23}(\Delta t)^{\alpha}}{n^{s(\alpha)}}\cdot 2\int\limits_0^{(\Delta t)^{-\alpha}} (1+r(\Delta t)^\alpha)dr =\frac{3C_{23}}{n^{s(\alpha)}}.
$$
For the second term in (5.4), we use the estimate
$$
|v_n|\leq\frac{C_{26}}{|\lambda|(\Delta t)^\alpha n^{s(\alpha)}},
$$
which follows by Lemmas 3.1 and 3.2. We obtain
$$
\int\limits_{\Gamma_2} \frac{|v_n(\lambda)|}{1+|\lambda|}|d\lambda| \leq\frac{C_{26}}{(\Delta t)^\alpha n^{s(\alpha)}} \int\limits_{\Gamma_2} \frac{|d\lambda|}{|\lambda|(1+|\lambda|)} \leq\frac{2C_{26}}{(\Delta t)^\alpha n^{s(\alpha)}} \int\limits_{(\Delta t)^{-\alpha}}^{+\infty} \frac{dr}{r^2} =\frac{2C_{26}}{n^{s(\alpha)}}.
$$
For the third term in (5.4), from (2.13) we get
$$
|E_\alpha(\lambda t_n^\alpha)| \leq\frac{C_2}{|\lambda|t_n^\alpha} =\frac{C_2}{|\lambda| (\Delta t) ^{\alpha} n^\alpha};
$$
$$
\int\limits_{\Gamma_2} \frac{|E_\alpha(\lambda t_n^\alpha)|}{1+|\lambda|}|d\lambda| \leq\frac{C_2}{(\Delta t)^\alpha n^\alpha} \int\limits_{\Gamma_2}\frac{|d\lambda|}{|\lambda|(1+|\lambda|)} \leq\frac{2C_2}{(\Delta t)^\alpha n^\alpha} \int\limits_{(\Delta t)^{-\alpha}}^{+\infty} \frac{dr}{r^2} =\frac{2C_2}{n^\alpha}.
$$
Substituting these estimates into (5.4), we conclude that
$$
\|u_n-u(t_n)\|\leq \frac{C_{27}\|f\|}{n^{s(\alpha)}}. \eqno{(5.5)}
$$

From (5.5) we get the pointwise convergence of the difference scheme (5.1) to the solution of (1.1). Let us explain this in more detail. Let $u_n=u_n^N$ be the approximation for $u(t_n)=u(nT/N)$ generated by (5.1), when the number of discretization steps is $N$. To approximate the values $u(t)$ at any given point $t\in[0,T]$, one can use $u_{[t/\Delta t]}^N$ with a sufficiently large $N$. Here, $[x]$ denotes the integer part of $x\in\mathbb{R}$. We claim that
$$
\forall t\in[0,T],\quad u_{[t/\Delta t]}^N \to u(t) \quad (N\to\infty). \eqno{(5.6)}
$$
At $t=0$, we have $u_{[t/\Delta t]}^N=u_0^N=f=u(0)$, hence (5.6) holds. Further assume that $t\in(0,T]$. Since $[t/\Delta t]\Delta t\to t$ when $N\to\infty$, and the function $u(t)$, $t\in[0,T]$ is continuous, then $u([t/\Delta t]\Delta t)\to u(t)$, as $N\to\infty$. Therefore, it is sufficient to verify that $\|u_{[t/\Delta t]}^N -u([t/\Delta t]\Delta t)\|\to 0$ when $N\to\infty$. But $\|u_{[t/\Delta t]}^N -u([t/\Delta t]\Delta t)\|$ is the value $\|u_n-u(t_n)\|$ in (5.5) with the choice $n=\Bigl[\frac{t}{\Delta t} \Bigr] =\Bigl[\frac{t}{T}N\Bigr]$. According to (5.5),
$$
\|u_{[t/\Delta t]}^N -u([t/\Delta t]\Delta t)\|\leq \frac{C_{27}\|f\|}{[t/\Delta t]^{s(\alpha)}} \leq C_{28}\|f\|(\Delta t)^{s(\alpha)}\to 0\quad (N\to\infty).
$$
Thus, we have established the pointwise convergence of the difference scheme (5.1) for solving (1.1) and proved the rate--of--convergence estimate
$$
\|u_{[t/\Delta t]}^N -u([t/\Delta t]\Delta t)\|\leq C_{28}\|f\|(\Delta t)^{s(\alpha)},\quad C_{28}=C_{28}(t). \eqno{(5.7)}
$$

For practical purposes it is important to ensure the stability of the scheme (5.1), i.e., to prove that the norms $\|u_n^N\|$ are uniformly bounded for all $N$ and $0\leq n\leq N$ by a value linear in $\|f\|$. We note that (5.5) implies $\|u_n-u(t_n)\|\equiv\|u_n^N-u(nT/N)\|\leq C_{27}\|f\|$ with a constant $C_{27}$ which does not depend on $N$, $n$, $\|f\|$. From Lemma 1.1 it follows that there exists a constant $C_{29}$ such that $\|u(t_n)\|\leq C_{29}\|f\|$. This means the stability of the method (5.1).

We have proved the following theorem which is the main result of the paper.

{\bf Theorem 5.1.} {\it Let Condition 1.1 be fulfilled. The finite difference scheme (5.1) is stable, and pointwisely converging in the sense of (5.6). For its convergence rate, the estimate (5.7) is true. In particular, the elements $u_N=u_N^N$ converge to $u(T)$ when $N\to\infty$. Moreover,
$$
\|u_N -u(T)\|\leq C_{28}\|f\|(\Delta t)^{s(\alpha)},
$$
where $s(\alpha)$ is defined in (3.17).}

\begin{center}
\textbf{6. Numerical experiments}
\end{center}

In this section we present results of numerical experiments with the difference scheme (5.1), as applied to the initial boundary value problem
\begin{equation}
\begin{split}
&\partial^\alpha_t x(s,t)=a_0^2\frac{\partial^2 x}{\partial
s^2}-b(s)\frac{\partial x}{\partial s}-c(s)x,\quad x=x(s,t); \\
&x(0,t)=0,\quad x(1,t)=0,\quad 0\leq t\leq T; \\
&x(s,0)=f(s),\quad 0\leq s\leq 1,\quad f(0)=f(1)=0.
\end{split}
\tag{6.1}
\end{equation}
Problem (6.1) has the form (1.1) with $u(t)=x(\cdot,t)$ and the operator $A:L_2[0,1]\to L_2[0,1]$ defined by
\begin{equation}
\begin{split}
[Au](s)&=a_0^2u^{\prime\prime}(s)-b(s)u^\prime(s)-c(s)u(s),\quad
s\in[0,1],\\
D(A)&=\{u\in H^2[0,1]\,|\, u(0)=u(1)=0\}.
\end{split}
\notag
\end{equation}
Here, the real functions $b\in C^1[0,1]$, $c\in C[0,1]$ are such that $A$
satisfies Condition 1.1, see Examples 6.1 and 6.2.

In accordance with (5.3), at each step we have to solve the following boundary value problem for the function $u_n(s)$, $s\in[0,1]$:
$$
(\Delta t)^\alpha\Bigl( a_0^2u_n^{\prime\prime}(s)-b(s)u_n^\prime(s)-c(s)u_n(s) \Bigr) - a_{nn}u_n(s) =\sum\limits_{j=0}^{n-1} a_{jn}u_j(s),\quad u_n(0)=u_n(1)=0. \eqno{(6.2)}
$$
The calculations were performed in Maple 15. At each iteration of (5.3), the boundary value problem (6.2) was solved on a grid on $[0,1]$ with the uniform partitioning in $2048$ intervals. Then we applied the cubic spline interpolation with the values $u_n(s)$ at $s=j/256$, $j=0,1,\dots,256$. The resulting function was used in the next iterations.

Along with (5.1), to solve the problem (6.1) we have also used the iterative method from [1,p.30]:
\begin{equation}
(\Delta t)^{-\alpha}\sum\limits_{j=0}^n (-1)^j
\frac{\alpha(\alpha-1)\cdot\dots\cdot(\alpha-j+1)}{j!}
(u_{n-j}-f)=Au_n,\quad 1\leq n\leq N;\quad u_0=f.
\tag{6.3}
\end{equation}
The finite difference scheme (6.3), as well as (5.1), is applicable to a wide class of problems (1.1), in particular to the problems (6.1). The implementation of (6.3) involves solution of boundary value problems analogous to (6.2) at each step.

{\bf Example 6.1.} We put $a_0=1$, $b(s)\equiv 0$, $c(s)\equiv 0$ in (6.1). The corresponding operator $A$ is self-adjoint and therefore satisfies Condition 1.1 for any $\varphi_0\in(0,\pi/2)$. We consider the problem of finding $x(s,T)$, $s\in[0,1]$, given the function $f(s)$, $s\in[0,1]$. Exact solution of the problem is easily found by the separation of variables:
$$
x(s,t)=\sum\limits_{n=1}^\infty\mathcal{S}_n(s)\mathcal{T}_n(t),\quad \mathcal{S}_n(s)=\sin\pi ns,\quad \mathcal{T}_n(t)=\mathcal{T}_n(0)E_\alpha(-\pi^2 n^2 t^\alpha),\quad n=1,2,\dots;
$$
$$
x(s,T)=\sum\limits_{n=1}^\infty \Bigl( f_n E_\alpha(-\pi^2 n^2 T^\alpha) \Bigr) \sin\pi ns,\quad f_n=2\int\limits_0^1 f(s)\sin\pi ns ds. \eqno{(6.4)}
$$
Explicit formulas for the solutions of (6.1) with an arbitrary self-adjoint differential operator of second order in the right part is given in [22].

In the numerical experiments we put $T=1$ while the function $f$, the value $\alpha$ and the number of iterations $N$ varied. The exact solution $u(T)=x(\cdot,T)$ of (6.1) is calculated according to (6.4). It is compared with the approximate solution $u_N$ obtained by the scheme (5.1) and with the approximate solution $\widetilde u_N$ obtained by (6.3) with the same $N$, the number of iterations. Numerical results are given in Table 1. In each test, we calculate the absolute error $\|u_N-u(T)\|$ of the scheme (5.1), the relative error $\frac{\|u_N-u(T)\|}{\|u(T)\|}$, the error $\|\widetilde u_N-u(T)\|$ of the method (6.3) and the corresponding relative error $\frac{\|\widetilde u_N-u(T)\|}{\|u(T)\|}$. All the norms are taken in $L_2[0,1]$. We also indicate the distance $\|u(T)-f\|$ in $L_2[0,1]$ between the initial element $f$ and the exact solution $u(T)$. As a rule, it is much more than the errors of both methods.

\begin{table}
\caption{Numerical results in Example 6.1}
\begin{center}
\begin{tabular}{|c|c|c|c|c|c|c|c|c|}\hline
No. & $f(s)$ & $\alpha$ & $N$ & $\|u_N-u(T)\|$ & $\frac{\|u_N-u(T)\|}{\|u(T)\|}$ & $\|\widetilde u_N-u(T)\|$ & $\frac{\|\widetilde u_N-u(T)\|}{\|u(T)\|}$ & $\|u(T)-f\|$\\ \hline
1 & $s^2(s-1)$ & $0.25$ & $5$ & $6.46\cdot 10^{-5}$ & $9.14\cdot 10^{-3}$ & $2.13\cdot 10^{-4}$ & $3.02\cdot 10^{-2}$ & $9.08\cdot 10^{-2}$\\ \hline
2 & --- & --- & $100$ & $3.34\cdot 10^{-6}$ & $4.72\cdot 10^{-4}$ & $1.06\cdot 10^{-5}$ & $1.50\cdot 10^{-3}$ & --- \\ \hline
3 & --- & $0.75$ & $5$ & $2.86\cdot 10^{-4}$ & $1.01\cdot 10^{-1}$ & $5.58\cdot 10^{-4}$ & $1.96\cdot 10^{-1}$ & $9.49\cdot 10^{-2}$ \\ \hline
4 & --- & --- & $100$ & $7.98\cdot 10^{-6}$ & $2.80\cdot 10^{-3}$ & $2.35\cdot 10^{-5}$ & $8.24\cdot 10^{-3}$ & --- \\ \hline
5 & $\sin 2\pi s$ & $0.25$ & $5$ & $1.48\cdot 10^{-4}$ & $1.03\cdot 10^{-2}$ & $4.68\cdot 10^{-4}$ & $3.26\cdot 10^{-2}$ & $6.93\cdot 10^{-1}$ \\ \hline
6 & --- & --- & $100$ & $6.67\cdot 10^{-6}$ & $4.64\cdot 10^{-4}$ & $2.20\cdot 10^{-5}$ & $1.53\cdot 10^{-3}$ & --- \\ \hline
7 & --- & $0.75$ & $5$ & $3.94\cdot 10^{-4}$ & $7.76\cdot 10^{-2}$ & $8.37\cdot 10^{-4}$ & $1.65\cdot 10^{-1}$ & $7.02\cdot 10^{-1}$ \\ \hline
8 & --- & --- & $100$ & $1.40\cdot 10^{-5}$ & $2.77\cdot 10^{-3}$ & $3.52\cdot 10^{-5}$ & $6.94\cdot 10^{-3}$ & --- \\ \hline
\end{tabular}
\end{center}
\end{table}

\begin{table}
\caption{Numerical results in Example 6.2}
\begin{center}
\begin{tabular}{|c|c|c|c|c|c|}\hline
No. & $f(s)$ & $\alpha$ & $N$ & $M$ & $\|u_N-\widetilde u_M\|$ \\ \hline
1 & $s^2(s-1)$ & $0.25$ & $5$ & $5$ & $4.00\cdot 10^{-4}$ \\ \hline
2 & --- & --- & $5$ & $100$ & $3.25\cdot 10^{-5}$ \\ \hline
3 & --- & --- & $100$ & $100$ & $2.10\cdot 10^{-5}$ \\ \hline
4 & --- & $0.75$ & $5$ & $5$ & $5.25\cdot 10^{-4}$ \\ \hline
5 & --- & --- & $5$ & $100$ & $4.51\cdot 10^{-4}$ \\ \hline
6 & --- & --- & $100$ & $100$ & $4.00\cdot 10^{-5}$ \\ \hline
7 & $\sin 2\pi s$ & $0.25$ & $5$ & $5$ & $3.55\cdot 10^{-3}$ \\ \hline
8 & --- & --- & $5$ & $100$ & $4.53\cdot 10^{-4}$ \\ \hline
9 & --- & --- & $100$ & $100$ & $1.89\cdot 10^{-4}$ \\ \hline
10 & --- & $0.75$ & $5$ & $5$ & $5.28\cdot 10^{-3}$ \\ \hline
11 & --- & --- & $5$ & $100$ & $6.11\cdot 10^{-3}$ \\ \hline
12 & --- & --- & $100$ & $100$ & $4.46\cdot 10^{-4}$ \\ \hline
\end{tabular}
\end{center}
\end{table}

In all our experiments, the method (5.1) has shown slightly better results than (6.3), although the calculation times were close. The absolute advantage of the scheme (5.1) is the presence of the theoretical error estimate, see Theorem 5.1, which guarantees the efficiency of (5.1), as applied to other similar problems.

{\bf Example 6.2.} Let us now consider the problem (6.1) with $a_0=0.1$, $b(s)=0.02 s$, $c(s)=s(1-s)+0.02$. Functions $b(s)$ and $c(s)$ satisfy the condition $c(s)-1/2 b^\prime(s)\geq\epsilon>0$ $\forall s\in[0,1]$. According to [23,Ch.5; 24], Condition 1.1 is fulfilled with $\varphi_0\in(\varphi_0^*,\pi/2)$, where
$$
\tan\varphi_0^*=\max\limits_{s\in[0,1]}\biggl(|b(s)|\max\biggl\{\frac{1}{2a_0^2},\frac{1}{2c(s)-b^\prime(s)}
\biggr\}\biggr).
$$

As above, we let $T=1$. In contrast to Example 6.1, for the solution $u(T)$ there is no explicit formula. Therefore, we can compare results of (5.1) with those obtained by (6.3). In our experiments, the function $f$, the value $\alpha$, the number of iterations $N$ of the method (5.1) and the number of iterations $M$ of the method (6.3) are varied, see Table 2. We give the distance in $L_2[0,1]$ between the approximate solutions $u_N$ and $\widetilde u_M$ obtained by (5.1) and (6.2).

Our calculations indicate that $\|u_N-\widetilde u_M\|\to 0$, as $N,M\to\infty$. Together with the data from Example 6.1, the results confirm the practical applicability of the difference scheme (5.1) studied in this article.

\medskip

\centerline{{\bf References}}

1. {\it Bajlekova E.G.} Fractional Evolution Equations in Banach
Spaces. Pleven: Eindhoven University of Technology, 2001.

2. {\it Kilbas A.A., Srivastava H.M., Trujillo J.J.} Theory and
Applications of Fractional Differential Equations. Amsterdam: Elsevier, 2006.

3. {\it Kochubey A.N.} A Cauchy problem for evolution equations of fractional order. {\it Differential Equations.} 1989. V.25. No.8. P.967--974.

4. {\it Taukenova F.I., Shkhanukov--Lafishev M.Kh.} Difference methods for solving boundary value problems for fractional differential equations. {\it Computational Mathematics and Mathematical Physics.} 2006. V.46. No.10. P.1785--1795.

5. {\it Li C., Zeng F.} The finite difference methods for fractional ordinary differential equations. {\it Numerical Functional Analysis and Optimization.} 2013. V.34. No.2. P.149--179.

6. {\it Abrashina--Zhadaeva N.G., Timoshchenko I.A.} Finite--difference schemes for a diffusion equation with fractional derivatives in a multidimensional domain. {\it Differential Equations.} 2013. V.49. No.7. P.789--795.

7. {\it Lafisheva M.M., Shkhanukov--Lafishev M.Kh.} Locally one--dimensional difference schemes for the fractional order diffusion equation. {\it Computational Mathematics and Mathematical Physics.} 2008. V.48. No.10. P.1875--1884.

8. {\it Alikhanov A.A.} Stability and convergence of difference schemes for boundary value problems for the fractional--order diffusion equation. {\it Computational Mathematics and Mathematical Physics.} 2016. V.56. No.4. P.561--575.

9. {\it Liu R., Li M., Pastor J., Piskarev S.I.} On the approximation of fractional resolution families. {\it Differential Equations.} 2014. V.50. No.7. P.927--937.

10. {\it Haase M.} The Functional Calculus for Sectorial Operators.
Basel: Birkh\"{a}user, 2006.

11. {\it Kokurin M.M.} The uniqueness of a solution to the inverse Cauchy problem for a fractional differential equation in a Banach space. {\it Russian Mathematics.} 2013. V.57. No.12. P.16--30.

12. {\it Kokurin M.M.} On the optimization of the rate--of--convergence estimates for some classes of difference schemes for solving ill--posed Cauchy problems. {\it Computational methods and programming.} 2013. V.14. P.58--76. (In Russian)

13. {\it Kokurin M.M.} Difference schemes for solving the Cauchy problem for a second--order operator differential equation. {\it Computational Mathematics and Mathematical Physics.} 2014. V.54. No.4. P.582--597.

14. {\it Lubich Ch.} Discretized fractional calculus. {\it SIAM Journal on Mathematical Analysis.} 1986. V.17. No.3. P.704--719.

15. {\it Dzhrbashyan M.M.} Integral Transforms and Representations of Functions in the Complex Domain. Moscow: Nauka, 1966. (In Russian)

16. {\it Popov A.Yu., Sedletskii A.M.} Distribution of roots of Mittag--Leffler functions. {\it Journal of Mathematical Sciences.} 2013. V.190. No.2. P.209--409.

17. {\it Vainberg M.M.} Variational Methods and Method of Monotone Operators. New York: Wiley, 1973.

18. {\it Trenogin V.A.} Functional Analysis. Moscow: FIZMATLIT, 2007. (In Russian)

19. {\it Radzievskaya E.I., Radzievskii G.V.} The remainder term of the Taylor expansion for a holomorphic function is representable in Lagrange form. {\it Siberian Mathematical Journal.} 2003. V.44. No.2. P.322--331.

20. {\it Bateman G., Erdelyi A.} Higher Transcendental Functions. V.1. New York, Toronto, London: Graw--Hill, 1953.

21. {\it Prudnikov A.P., Brychkov Y.A., Marichev O.I.} Integrals and Series. V.3. Special Functions. Additional Chapters. Moscow: FIZMATLIT, 2003. (In Russian)

22. {\it Sakamoto K., Yamamoto M.} Initial value/boundary value problems for fractional diffusion--wave equations and applications to some inverse problems. {\it Journal of Mathematical Analysis and Applications.} 2011. V.382. P.426--447.

23. {\it Kato T.} Perturbation Theory for Linear Operators. Berlin, Heidelberg, New York: Springer--Verlag, 1966.

24. {\it Bakushinskii A.B., Kokurin M.M. Kokurin M.Yu.} On a class of finite--difference schemes for solving ill--posed Cauchy problems in Banach spaces. {\it Computational Mathematics and Mathematical Physics.} 2012. V.52. No.3. P.411--426.

\end{document}